\documentclass{gtart}

\def\ifplaintex{\expandafter\ifx\csname documentclass\endcsname\relax}

\def\gtp{{\mathsurround=0pt\it $\cal G\mskip-2mu$eometry \&\ 
$\cal T\!\!$opology $\cal P\!$ublications}}  

\def\recd{{\small Received:\qua\receiveddate\ifx\reviseddate\relax
\else\qquad Revised:\qua\reviseddate\fi\par}} 

\def\lognumber#1{\def\thelognumber{#1}}
\def\volumenumber#1{\def\thevolumenumber{#1}}
\def\volumeyear#1{\def\thevolumeyear{#1}}
\def\papernumber#1{\def\thepapernumber{#1}}
\def\pagenumbers#1#2{\def\startpage{#1}\def\finishpage{#2}}
\def\published#1{\def\publishdate{#1}}

\def\received#1{\def\receiveddate{#1}}

\def\accepted#1{\def\accepteddate{#1}}

\long\def\asciiabstract#1{\long\def\theasciiabstract{#1}}

\let\\\par\let\thelognumber\relax\let\thevolumenumber\relax
\let\thepapernumber\relax\let\thevolumeyear\relax\let\startpage\relax
\let\finishpage\relax\let\publishdate\relax\let\receiveddate\relax
\let\reviseddate\relax\let\accepteddate\relax\let\theasciititle\relax
\let\theasciiauthors\relax
\let\theasciiabstract\relax

\let\theasciiemail\relax

\ifplaintex
\font\logobig=cmssbx10 scaled 3836
\font\logomed=cmssbx10 scaled 2557
\else
\font\logobig=cmssbx10 scaled 4200
\font\logomed=cmssbx10 scaled 2800
\fi

\long\def\makeagttitle{   
\count0=\startpage
\agt\hfill      
\hbox to 45truept{\vbox to 0pt{\vglue -13truept{\logomed A\kern -.37em{\logobig 
T}\kern -.38em G}\vss}\hss}
\break
{\small Volume \thevolumenumber\ (\thevolumeyear)
\startpage--\finishpage\nl
Published: \publishdate}

\vglue .25truein

{\parskip=0pt\leftskip 0pt plus
1fil\def\\{\par\smallskip}{\Large\bf\thetitle}\par\medskip} \vglue
0.05truein

{\parskip=0pt\leftskip 0pt plus 1fil\def\\{\par}{\sc\theauthors}
\par\medskip}%
 
\vglue 0.03truein 

{\small\leftskip 25truept\rightskip 25truept{\bf Abstract}\stdspace\theabstract

{\bf AMS Classification}\stdspace\theprimaryclass
\ifx\thesecondaryclass\relax\else; \thesecondaryclass\fi\par
{\bf Keywords}\stdspace \thekeywords\par}\vglue 7truept

}   

\ifplaintex
\font\phead=cmsl9 scaled 950
\font\pnum=cmbx10 scaled 913
\font\pfoot=cmsl9 scaled 950
\headline{\vbox to 0pt{\vskip -4.5mm\line{\small\phead\ifnum
\count0=\startpage ISSN 1472-2739 (on-line) 1472-2747 (printed)
\hfill {\pnum\folio}\else\ifodd\count0\def\\{ }%
\ifx\theshorttitle\relax\thetitle\else\theshorttitle\fi\hfill{\pnum\folio}
\else\def\\{ and }{\pnum\folio}\hfill\ifx\theshortauthors\relax\theauthors
\else\theshortauthors\fi\fi\fi}\vss}}
\footline{\vbox to 0pt{\vglue 0mm\line{\small\pfoot\ifnum\count0=\startpage
\copyright\ \gtp\hfill\else
\agt, Volume \thevolumenumber\ (\thevolumeyear)\hfill\fi}\vss}}
\else
\font=cmsl9 scaled 1050
\font\lnum=cmbx10 
\font=cmsl9 scaled 1050
\makeatletter
\def\@oddhead{{\small\lhead\ifnum\count0=\startpage ISSN 1472-2739 
(on-line) 1472-2747 (printed)\hfill {\lnum\number\count0}\else\ifodd\count0
\def\\{ }\ifx\theshorttitle\relax \thetitle \else\theshorttitle\fi\hfill
{\lnum\number\count0}\else\def\\{ and }{\lnum\number\count0}
\hfill\ifx\theshortauthors\relax 
\theauthors\else\theshortauthors\fi\fi\fi}}\def\@evenhead{\@oddhead}
\def\@oddfoot{\small\lfoot\ifnum\count0=\startpage\copyright\ \gtp\hfill\else
\agt, Volume \thevolumenumber\ (\thevolumeyear)\hfill\fi}
\def\@evenfoot{\@oddfoot}
\makeatother
\fi
\let\maketitlepage\makeagttitle
\let\makeshorttitle\maketitlepage
\let\maketitle\maketitlepage

\newwrite\gtoutfile
\long\gdef\makeheadfile{  
{\def\\{, }\def\s{ }
\immediate\openout\gtoutfile head.xxx
\immediate\write\gtoutfile{To: math@arxiv.org}
\immediate\write\gtoutfile{Subject: put OR rep NNNNN:ppppp}
\immediate\write\gtoutfile{--text follows this line--}
\immediate\write\gtoutfile{Proxy-for: \ifx\theasciiauthors\relax
\theauthors\else\theasciiauthors\fi\s<\ifx\theasciiemail\relax\theemail\else\theasciiemail\fi>}
\immediate\write\gtoutfile{\noexpand\\}
\immediate\write\gtoutfile{Authors: \ifx\theasciiauthors\relax
\theauthors\else\theasciiauthors\fi}
{\def\\{ }\immediate\write\gtoutfile{Title: \ifx\theasciititle\relax
\thetitle\else\theasciititle\fi}}
\immediate\write\gtoutfile{Subj-class: GT or SG, GR etc}
\immediate\write\gtoutfile{MSC-class: \theprimaryclass\ifx\thesecondaryclass\relax\else, \thesecondaryclass\fi}
\immediate\write\gtoutfile{Journal-ref: Algebr. Geom. Topol. \thevolumenumber\s
(\thevolumeyear) \startpage-\finishpage}
\immediate\write\gtoutfile{Comments: Published by Algebraic and
Geometric Topology at}
\immediate\write\gtoutfile{\s\s\s  http://www.maths.warwick.ac.uk/agt/AGTVol\thevolumenumber/agt-\thevolumenumber-\thepapernumber.abs.html}
\immediate\write\gtoutfile{\noexpand\\}
\immediate\write\gtoutfile{}
\ifx\theasciiabstract\relax
\immediate\write\gtoutfile{\theabstract}\else
\immediate\write\gtoutfile{\theasciiabstract}\fi
\immediate\write\gtoutfile{}
\immediate\write\gtoutfile{\noexpand\\}
\immediate\write\gtoutfile{}
\immediate\closeout\gtoutfile}}  

\def\maketitlepage{\makeagttitle\makeheadfile}
\let\makeshorttitle\maketitlepage
\let\maketitle\maketitlepage

\def\ifplaintex{\expandafter\ifx\csname documentclass\endcsname\relax}

\def\gtp{{\mathsurround=0pt\it $\cal G\mskip-2mu$eometry \&\ 
$\cal T\!\!$opology $\cal P\!$ublications}}  

\def\recd{{\small Received:\qua\receiveddate\ifx\reviseddate\relax
\else\qquad Revised:\qua\reviseddate\fi\par}} 

\def\lognumber#1{\def\thelognumber{#1}}
\def\volumenumber#1{\def\thevolumenumber{#1}}
\def\volumeyear#1{\def\thevolumeyear{#1}}
\def\papernumber#1{\def\thepapernumber{#1}}
\def\pagenumbers#1#2{\def\startpage{#1}\def\finishpage{#2}}
\def\published#1{\def\publishdate{#1}}

\def\received#1{\def\receiveddate{#1}}

\def\accepted#1{\def\accepteddate{#1}}

\long\def\asciiabstract#1{\long\def\theasciiabstract{#1}}

\let\\\par\let\thelognumber\relax\let\thevolumenumber\relax
\let\thepapernumber\relax\let\thevolumeyear\relax\let\startpage\relax
\let\finishpage\relax\let\publishdate\relax\let\receiveddate\relax
\let\reviseddate\relax\let\accepteddate\relax\let\theasciititle\relax
\let\theasciiauthors\relax
\let\theasciiabstract\relax

\let\theasciiemail\relax

\ifplaintex
\font\logobig=cmssbx10 scaled 3836
\font\logomed=cmssbx10 scaled 2557
\else
\font\logobig=cmssbx10 scaled 4200
\font\logomed=cmssbx10 scaled 2800
\fi

\long\def\makeagttitle{   
\count0=\startpage
\agt\hfill      
\hbox to 45truept{\vbox to 0pt{\vglue -13truept{\logomed A\kern -.37em{\logobig 
T}\kern -.38em G}\vss}\hss}
\break
{\small Volume \thevolumenumber\ (\thevolumeyear)
\startpage--\finishpage\nl
Published: \publishdate}

\vglue .25truein

{\parskip=0pt\leftskip 0pt plus
1fil\def\\{\par\smallskip}{\Large\bf\thetitle}\par\medskip} \vglue
0.05truein

{\parskip=0pt\leftskip 0pt plus 1fil\def\\{\par}{\sc\theauthors}
\par\medskip}%
 
\vglue 0.03truein 

{\small\leftskip 25truept\rightskip 25truept{\bf Abstract}\stdspace\theabstract

{\bf AMS Classification}\stdspace\theprimaryclass
\ifx\thesecondaryclass\relax\else; \thesecondaryclass\fi\par
{\bf Keywords}\stdspace \thekeywords\par}\vglue 7truept

}   

\ifplaintex
\font\phead=cmsl9 scaled 950
\font\pnum=cmbx10 scaled 913
\font\pfoot=cmsl9 scaled 950
\headline{\vbox to 0pt{\vskip -4.5mm\line{\small\phead\ifnum
\count0=\startpage ISSN 1472-2739 (on-line) 1472-2747 (printed)
\hfill {\pnum\folio}\else\ifodd\count0\def\\{ }%
\ifx\theshorttitle\relax\thetitle\else\theshorttitle\fi\hfill{\pnum\folio}
\else\def\\{ and }{\pnum\folio}\hfill\ifx\theshortauthors\relax\theauthors
\else\theshortauthors\fi\fi\fi}\vss}}
\footline{\vbox to 0pt{\vglue 0mm\line{\small\pfoot\ifnum\count0=\startpage
\copyright\ \gtp\hfill\else
\agt, Volume \thevolumenumber\ (\thevolumeyear)\hfill\fi}\vss}}
\else
\font=cmsl9 scaled 1050
\font\lnum=cmbx10 
\font=cmsl9 scaled 1050
\makeatletter
\def\@oddhead{{\small\lhead\ifnum\count0=\startpage ISSN 1472-2739 
(on-line) 1472-2747 (printed)\hfill {\lnum\number\count0}\else\ifodd\count0
\def\\{ }\ifx\theshorttitle\relax \thetitle \else\theshorttitle\fi\hfill
{\lnum\number\count0}\else\def\\{ and }{\lnum\number\count0}
\hfill\ifx\theshortauthors\relax 
\theauthors\else\theshortauthors\fi\fi\fi}}\def\@evenhead{\@oddhead}
\def\@oddfoot{\small\lfoot\ifnum\count0=\startpage\copyright\ \gtp\hfill\else
\agt, Volume \thevolumenumber\ (\thevolumeyear)\hfill\fi}
\def\@evenfoot{\@oddfoot}
\makeatother
\fi
\let\maketitlepage\makeagttitle
\let\makeshorttitle\maketitlepage
\let\maketitle\maketitlepage

\newwrite\gtoutfile
\long\gdef\makeheadfile{  
{\def\\{, }\def\s{ }
\immediate\openout\gtoutfile head.xxx
\immediate\write\gtoutfile{To: math@arxiv.org}
\immediate\write\gtoutfile{Subject: put OR rep NNNNN:ppppp}
\immediate\write\gtoutfile{--text follows this line--}
\immediate\write\gtoutfile{Proxy-for: \ifx\theasciiauthors\relax
\theauthors\else\theasciiauthors\fi\s<\ifx\theasciiemail\relax\theemail\else\theasciiemail\fi>}
\immediate\write\gtoutfile{\noexpand\\}
\immediate\write\gtoutfile{Authors: \ifx\theasciiauthors\relax
\theauthors\else\theasciiauthors\fi}
{\def\\{ }\immediate\write\gtoutfile{Title: \ifx\theasciititle\relax
\thetitle\else\theasciititle\fi}}
\immediate\write\gtoutfile{Subj-class: GT or SG, GR etc}
\immediate\write\gtoutfile{MSC-class: \theprimaryclass\ifx\thesecondaryclass\relax\else, \thesecondaryclass\fi}
\immediate\write\gtoutfile{Journal-ref: Algebr. Geom. Topol. \thevolumenumber\s
(\thevolumeyear) \startpage-\finishpage}
\immediate\write\gtoutfile{Comments: Published by Algebraic and
Geometric Topology at}
\immediate\write\gtoutfile{\s\s\s  http://www.maths.warwick.ac.uk/agt/AGTVol\thevolumenumber/agt-\thevolumenumber-\thepapernumber.abs.html}
\immediate\write\gtoutfile{\noexpand\\}
\immediate\write\gtoutfile{}
\ifx\theasciiabstract\relax
\immediate\write\gtoutfile{\theabstract}\else
\immediate\write\gtoutfile{\theasciiabstract}\fi
\immediate\write\gtoutfile{}
\immediate\write\gtoutfile{\noexpand\\}
\immediate\write\gtoutfile{}
\immediate\closeout\gtoutfile}}  

\def\maketitlepage{\makeagttitle\makeheadfile}
\let\makeshorttitle\maketitlepage
\let\maketitle\maketitlepage

\def\ifplaintex{\expandafter\ifx\csname documentclass\endcsname\relax}

\def\gtp{{\mathsurround=0pt\it $\cal G\mskip-2mu$eometry \&\ 
$\cal T\!\!$opology $\cal P\!$ublications}}  

\def\recd{{\small Received:\qua\receiveddate\ifx\reviseddate\relax
\else\qquad Revised:\qua\reviseddate\fi\par}} 

\def\lognumber#1{\def\thelognumber{#1}}
\def\volumenumber#1{\def\thevolumenumber{#1}}
\def\volumeyear#1{\def\thevolumeyear{#1}}
\def\papernumber#1{\def\thepapernumber{#1}}
\def\pagenumbers#1#2{\def\startpage{#1}\def\finishpage{#2}}
\def\published#1{\def\publishdate{#1}}

\def\received#1{\def\receiveddate{#1}}

\def\accepted#1{\def\accepteddate{#1}}

\long\def\asciiabstract#1{\long\def\theasciiabstract{#1}}

\let\\\par\let\thelognumber\relax\let\thevolumenumber\relax
\let\thepapernumber\relax\let\thevolumeyear\relax\let\startpage\relax
\let\finishpage\relax\let\publishdate\relax\let\receiveddate\relax
\let\reviseddate\relax\let\accepteddate\relax\let\theasciititle\relax
\let\theasciiauthors\relax
\let\theasciiabstract\relax

\let\theasciiemail\relax

\ifplaintex
\font\logobig=cmssbx10 scaled 3836
\font\logomed=cmssbx10 scaled 2557
\else
\font\logobig=cmssbx10 scaled 4200
\font\logomed=cmssbx10 scaled 2800
\fi

\long\def\makeagttitle{   
\count0=\startpage
\agt\hfill      
\hbox to 45truept{\vbox to 0pt{\vglue -13truept{\logomed A\kern -.37em{\logobig 
T}\kern -.38em G}\vss}\hss}
\break
{\small Volume \thevolumenumber\ (\thevolumeyear)
\startpage--\finishpage\nl
Published: \publishdate}

\vglue .25truein

{\parskip=0pt\leftskip 0pt plus
1fil\def\\{\par\smallskip}{\Large\bf\thetitle}\par\medskip} \vglue
0.05truein

{\parskip=0pt\leftskip 0pt plus 1fil\def\\{\par}{\sc\theauthors}
\par\medskip}%
 
\vglue 0.03truein 

{\small\leftskip 25truept\rightskip 25truept{\bf Abstract}\stdspace\theabstract

{\bf AMS Classification}\stdspace\theprimaryclass
\ifx\thesecondaryclass\relax\else; \thesecondaryclass\fi\par
{\bf Keywords}\stdspace \thekeywords\par}\vglue 7truept

}   

\ifplaintex
\font\phead=cmsl9 scaled 950
\font\pnum=cmbx10 scaled 913
\font\pfoot=cmsl9 scaled 950
\headline{\vbox to 0pt{\vskip -4.5mm\line{\small\phead\ifnum
\count0=\startpage ISSN 1472-2739 (on-line) 1472-2747 (printed)
\hfill {\pnum\folio}\else\ifodd\count0\def\\{ }%
\ifx\theshorttitle\relax\thetitle\else\theshorttitle\fi\hfill{\pnum\folio}
\else\def\\{ and }{\pnum\folio}\hfill\ifx\theshortauthors\relax\theauthors
\else\theshortauthors\fi\fi\fi}\vss}}
\footline{\vbox to 0pt{\vglue 0mm\line{\small\pfoot\ifnum\count0=\startpage
\copyright\ \gtp\hfill\else
\agt, Volume \thevolumenumber\ (\thevolumeyear)\hfill\fi}\vss}}
\else
\font=cmsl9 scaled 1050
\font\lnum=cmbx10 
\font=cmsl9 scaled 1050
\makeatletter
\def\@oddhead{{\small\lhead\ifnum\count0=\startpage ISSN 1472-2739 
(on-line) 1472-2747 (printed)\hfill {\lnum\number\count0}\else\ifodd\count0
\def\\{ }\ifx\theshorttitle\relax \thetitle \else\theshorttitle\fi\hfill
{\lnum\number\count0}\else\def\\{ and }{\lnum\number\count0}
\hfill\ifx\theshortauthors\relax 
\theauthors\else\theshortauthors\fi\fi\fi}}\def\@evenhead{\@oddhead}
\def\@oddfoot{\small\lfoot\ifnum\count0=\startpage\copyright\ \gtp\hfill\else
\agt, Volume \thevolumenumber\ (\thevolumeyear)\hfill\fi}
\def\@evenfoot{\@oddfoot}
\makeatother
\fi
\let\maketitlepage\makeagttitle
\let\makeshorttitle\maketitlepage
\let\maketitle\maketitlepage

\newwrite\gtoutfile
\long\gdef\makeheadfile{  
{\def\\{, }\def\s{ }
\immediate\openout\gtoutfile head.xxx
\immediate\write\gtoutfile{To: math@arxiv.org}
\immediate\write\gtoutfile{Subject: put OR rep NNNNN:ppppp}
\immediate\write\gtoutfile{--text follows this line--}
\immediate\write\gtoutfile{Proxy-for: \ifx\theasciiauthors\relax
\theauthors\else\theasciiauthors\fi\s<\ifx\theasciiemail\relax\theemail\else\theasciiemail\fi>}
\immediate\write\gtoutfile{\noexpand\\}
\immediate\write\gtoutfile{Authors: \ifx\theasciiauthors\relax
\theauthors\else\theasciiauthors\fi}
{\def\\{ }\immediate\write\gtoutfile{Title: \ifx\theasciititle\relax
\thetitle\else\theasciititle\fi}}
\immediate\write\gtoutfile{Subj-class: GT or SG, GR etc}
\immediate\write\gtoutfile{MSC-class: \theprimaryclass\ifx\thesecondaryclass\relax\else, \thesecondaryclass\fi}
\immediate\write\gtoutfile{Journal-ref: Algebr. Geom. Topol. \thevolumenumber\s
(\thevolumeyear) \startpage-\finishpage}
\immediate\write\gtoutfile{Comments: Published by Algebraic and
Geometric Topology at}
\immediate\write\gtoutfile{\s\s\s  http://www.maths.warwick.ac.uk/agt/AGTVol\thevolumenumber/agt-\thevolumenumber-\thepapernumber.abs.html}
\immediate\write\gtoutfile{\noexpand\\}
\immediate\write\gtoutfile{}
\ifx\theasciiabstract\relax
\immediate\write\gtoutfile{\theabstract}\else
\immediate\write\gtoutfile{\theasciiabstract}\fi
\immediate\write\gtoutfile{}
\immediate\write\gtoutfile{\noexpand\\}
\immediate\write\gtoutfile{}
\immediate\closeout\gtoutfile}}  

\def\maketitlepage{\makeagttitle\makeheadfile}
\let\makeshorttitle\maketitlepage
\let\maketitle\maketitlepage

\def\ifplaintex{\expandafter\ifx\csname documentclass\endcsname\relax}

\def\gtp{{\mathsurround=0pt\it $\cal G\mskip-2mu$eometry \&\ 
$\cal T\!\!$opology $\cal P\!$ublications}}  

\def\recd{{\small Received:\qua\receiveddate\ifx\reviseddate\relax
\else\qquad Revised:\qua\reviseddate\fi\par}} 

\def\lognumber#1{\def\thelognumber{#1}}
\def\volumenumber#1{\def\thevolumenumber{#1}}
\def\volumeyear#1{\def\thevolumeyear{#1}}
\def\papernumber#1{\def\thepapernumber{#1}}
\def\pagenumbers#1#2{\def\startpage{#1}\def\finishpage{#2}}
\def\published#1{\def\publishdate{#1}}

\def\received#1{\def\receiveddate{#1}}

\def\accepted#1{\def\accepteddate{#1}}

\long\def\asciiabstract#1{\long\def\theasciiabstract{#1}}

\let\\\par\let\thelognumber\relax\let\thevolumenumber\relax
\let\thepapernumber\relax\let\thevolumeyear\relax\let\startpage\relax
\let\finishpage\relax\let\publishdate\relax\let\receiveddate\relax
\let\reviseddate\relax\let\accepteddate\relax\let\theasciititle\relax
\let\theasciiauthors\relax
\let\theasciiabstract\relax

\let\theasciiemail\relax

\ifplaintex
\font\logobig=cmssbx10 scaled 3836
\font\logomed=cmssbx10 scaled 2557
\else
\font\logobig=cmssbx10 scaled 4200
\font\logomed=cmssbx10 scaled 2800
\fi

\long\def\makeagttitle{   
\count0=\startpage
\agt\hfill      
\hbox to 45truept{\vbox to 0pt{\vglue -13truept{\logomed A\kern -.37em{\logobig 
T}\kern -.38em G}\vss}\hss}
\break
{\small Volume \thevolumenumber\ (\thevolumeyear)
\startpage--\finishpage\nl
Published: \publishdate}

\vglue .25truein

{\parskip=0pt\leftskip 0pt plus
1fil\def\\{\par\smallskip}{\Large\bf\thetitle}\par\medskip} \vglue
0.05truein

{\parskip=0pt\leftskip 0pt plus 1fil\def\\{\par}{\sc\theauthors}
\par\medskip}%
 
\vglue 0.03truein 

{\small\leftskip 25truept\rightskip 25truept{\bf Abstract}\stdspace\theabstract

{\bf AMS Classification}\stdspace\theprimaryclass
\ifx\thesecondaryclass\relax\else; \thesecondaryclass\fi\par
{\bf Keywords}\stdspace \thekeywords\par}\vglue 7truept

}   

\ifplaintex
\font\phead=cmsl9 scaled 950
\font\pnum=cmbx10 scaled 913
\font\pfoot=cmsl9 scaled 950
\headline{\vbox to 0pt{\vskip -4.5mm\line{\small\phead\ifnum
\count0=\startpage ISSN 1472-2739 (on-line) 1472-2747 (printed)
\hfill {\pnum\folio}\else\ifodd\count0\def\\{ }%
\ifx\theshorttitle\relax\thetitle\else\theshorttitle\fi\hfill{\pnum\folio}
\else\def\\{ and }{\pnum\folio}\hfill\ifx\theshortauthors\relax\theauthors
\else\theshortauthors\fi\fi\fi}\vss}}
\footline{\vbox to 0pt{\vglue 0mm\line{\small\pfoot\ifnum\count0=\startpage
\copyright\ \gtp\hfill\else
\agt, Volume \thevolumenumber\ (\thevolumeyear)\hfill\fi}\vss}}
\else
\font=cmsl9 scaled 1050
\font\lnum=cmbx10 
\font=cmsl9 scaled 1050
\makeatletter
\def\@oddhead{{\small\lhead\ifnum\count0=\startpage ISSN 1472-2739 
(on-line) 1472-2747 (printed)\hfill {\lnum\number\count0}\else\ifodd\count0
\def\\{ }\ifx\theshorttitle\relax \thetitle \else\theshorttitle\fi\hfill
{\lnum\number\count0}\else\def\\{ and }{\lnum\number\count0}
\hfill\ifx\theshortauthors\relax 
\theauthors\else\theshortauthors\fi\fi\fi}}\def\@evenhead{\@oddhead}
\def\@oddfoot{\small\lfoot\ifnum\count0=\startpage\copyright\ \gtp\hfill\else
\agt, Volume \thevolumenumber\ (\thevolumeyear)\hfill\fi}
\def\@evenfoot{\@oddfoot}
\makeatother
\fi
\let\maketitlepage\makeagttitle
\let\makeshorttitle\maketitlepage
\let\maketitle\maketitlepage

\newwrite\gtoutfile
\long\gdef\makeheadfile{  
{\def\\{, }\def\s{ }
\immediate\openout\gtoutfile head.xxx
\immediate\write\gtoutfile{To: math@arxiv.org}
\immediate\write\gtoutfile{Subject: put OR rep NNNNN:ppppp}
\immediate\write\gtoutfile{--text follows this line--}
\immediate\write\gtoutfile{Proxy-for: \ifx\theasciiauthors\relax
\theauthors\else\theasciiauthors\fi\s<\ifx\theasciiemail\relax\theemail\else\theasciiemail\fi>}
\immediate\write\gtoutfile{\noexpand\\}
\immediate\write\gtoutfile{Authors: \ifx\theasciiauthors\relax
\theauthors\else\theasciiauthors\fi}
{\def\\{ }\immediate\write\gtoutfile{Title: \ifx\theasciititle\relax
\thetitle\else\theasciititle\fi}}
\immediate\write\gtoutfile{Subj-class: GT or SG, GR etc}
\immediate\write\gtoutfile{MSC-class: \theprimaryclass\ifx\thesecondaryclass\relax\else, \thesecondaryclass\fi}
\immediate\write\gtoutfile{Journal-ref: Algebr. Geom. Topol. \thevolumenumber\s
(\thevolumeyear) \startpage-\finishpage}
\immediate\write\gtoutfile{Comments: Published by Algebraic and
Geometric Topology at}
\immediate\write\gtoutfile{\s\s\s  http://www.maths.warwick.ac.uk/agt/AGTVol\thevolumenumber/agt-\thevolumenumber-\thepapernumber.abs.html}
\immediate\write\gtoutfile{\noexpand\\}
\immediate\write\gtoutfile{}
\ifx\theasciiabstract\relax
\immediate\write\gtoutfile{\theabstract}\else
\immediate\write\gtoutfile{\theasciiabstract}\fi
\immediate\write\gtoutfile{}
\immediate\write\gtoutfile{\noexpand\\}
\immediate\write\gtoutfile{}
\immediate\closeout\gtoutfile}}  

\def\maketitlepage{\makeagttitle\makeheadfile}
\let\makeshorttitle\maketitlepage
\let\maketitle\maketitlepage

\lognumber{12}
\volumenumber{1}
\volumeyear{2001}
\papernumber{12}
\pagenumbers{243}{270}
\received{14 November 2000}
\accepted{17 April 2001}
\published{23 April 2001}

\usepackage{amssymb,amsmath}
\usepackage{pb-diagram,lamsarrow,pb-lams,epsf}

\let\lbl\label
\theoremstyle{plain}

\newtheorem{theorem}{Theorem}
\newtheorem{proposition}{Proposition}[section]
\newtheorem{lemma}[proposition]{Lemma}
\newtheorem{corollary}[proposition]{Corollary}

\theoremstyle{definition}

\newtheorem{question}{Question}

\newtheorem{remark}[proposition]{Remark}

\newcommand{\psdraw}[2]{\centerline{\epsfxsize #2 \epsfbox{#1.eps}}}

\newlength{\standardunitlength}
\newcommand{\noi}{\noindent}

\renewcommand{\hom}{\operatorname{Hom}}
\renewcommand{\ker}{\operatorname{Ker}}
\newcommand{\im}{\operatorname{Im}}
\newcommand{\inte}{\operatorname{Int}}
\newcommand{\aut}{\operatorname{Aut}}
\newcommand{\ao}{\operatorname{Aut_0}}
\newcommand{\ai}{ \operatorname{Aut_1}}

\newcommand{\con}{\equiv}

\renewcommand{\rk}{\operatorname{rank}}
\newcommand{\auto}{\operatorname{Auto}}

\renewcommand{\sp}{\operatorname{Sp}}

\def\LL{\Lambda}

\def\SS{\Sigma}

\def\a{\alpha}
\def\b{\beta}

\def\th{\theta}

\def\t{\tau}
\newcommand{\s}{\sigma}
\def\g{\gamma}
\def\l{\lambda}

\def\jk{J^L_k}
\def\Lg{\mathcal L_g}
\def\Z{\mathbb Z}

\def\D{\mathsf D}

\def\C{\mathcal C}
\def\F{\mathcal F}

\def\A{\mathcal A}
\def\P{\mathcal P}

\def\L{{\mathsf L}}

\def\S{\mathcal S}

\def\G{\mathcal G}
\def\Q{\mathbb Q}

\def\M{\mathcal M}

\def\bd{\partial}
\def\Sg{\Sigma_{g,1}}
\def\sgo{\mathcal S_g^0}

\def\i{^{-1}}

\def\iso{\cong}

\def\sub{\subseteq}

\def\con{\equiv}
\def\pgf{\mathcal P_g^{\text{fr}}}
\def\sgf{\mathcal S_g^{\text{fr}}}
 
\def\Gg{\Gamma_{g,1}}
\def\gg{\Gamma_g}
\def\ggb{\Gamma_g^B}
\def\Hg{\mathcal H_g}
\def\bhg{\bar{\mathcal H_g}}
\def\pg{\mathcal P_g}
\def\sg{\mathcal S_g}
\def\ssg{\Sigma_g}
\def\dg{D_g}
\def\hp{\hat\Phi}

\def\Tg{\mathcal T\Hg}
\def\tg{T_g}
\def\hg{H_g}
\def\bg{\mathcal H_g^B}
\def\Bg{\Gg^B}
\def\cg{\mathcal C_g}

\begin{document}

\title
{Homology cylinders:  an enlargement of\\the mapping class group}
\author{Jerome Levine}
\address{Department of Mathematics, Brandeis University\\Waltham, 
MA 02454-9110, USA}        
\email{levine@brandeis.edu  }
\url{www.math.brandeis.edu/Faculty/jlevine/} 
\begin{abstract}
We consider a {\em homological} enlargement of the mapping class
group, defined by {\em homology cylinders} over a closed oriented
surface (up to homology cobordism). These are important model objects
in the recent Goussarov-Habiro theory of finite-type invariants of
$3$-manifolds. We study the structure of this group from several
directions: the relative weight filtration of Dennis Johnson, the
finite-type filtration of Goussarov-Habiro, and the relation to string
link concordance.

We also consider a new {\em Lagrangian} filtration of both the mapping class
group and the group of homology cylinders.
\end{abstract}
\asciiabstract{
We consider a homological enlargement of the mapping class
group, defined by homology cylinders over a closed oriented
surface (up to homology cobordism). These are important model objects
in the recent Goussarov-Habiro theory of finite-type invariants of
3-manifolds. We study the structure of this group from several
directions: the relative weight filtration of Dennis Johnson, the
finite-type filtration of Goussarov-Habiro, and the relation to string
link concordance.
We also consider a new Lagrangian filtration of both the mapping class
group and the group of homology cylinders.}

\primaryclass{57N10}
\secondaryclass{57M25}
\keywords{Homology cylinder, mapping class group, clasper, finite-type
invariant}

\makeshorttitle

\section{Introduction}

The mapping class group $\gg$ is the group 
of
diffeotopy
classes of orientation preserving diffeomorphisms of the closed oriented surface
$\ssg$
of genus $g$. There has been a great deal of work aimed at the determination of
the
algebraic structure of this group. For example, some years ago D. Johnson
defined a
filtration on $\gg$ (the
{\em relative weight filtration}) and observed that the associated graded group
is a Lie
subalgebra of a Lie algebra $\D(H)$ constructed explicitly from $H=H_1
(\ssg )$. Johnson, Morita and others (see \cite{J},\cite{Mo},\cite{M}) have
investigated
this Lie subalgebra but its precise determination is still open. R. Hain
\cite{Hn}
has studied the lower central series filtration of the Torelli group (the
subgroup of $\gg$ whose elements are homologically trivial) and found a simple
explicit presentation over $\Q$. We also mention the work of Oda (see 
\cite{O},\cite{L})
relating the pure braid group to $\gg$.

One reason for interest in $\gg$ is that it is related in an obvious way to the
structure of
$3$-manifolds via the Heegard construction. From this viewpoint,  the subgroup
$\ggb$
consisting of diffeomorphisms which extend over the handlebody $T_g$ of genus
$g$, and
the coset space $\gg /\ggb$, are of obvious interest.

In this note we propose a ``homological'' generalization of these groups, where
we replace
$\gg$ by a group $\Hg$ of homology bordism classes of homology cylinders over
$\ssg$. Homology cylinders have appeared and been studied in recent work of
Goussarov \cite{Go} and Habiro \cite{H}, as important model objects for their
new theory of finite-type invariants of general $3$-manifolds. We will see that
$\gg$ is a subgroup of $\Hg$. Furthermore the notion of
Heegard
construction translates to the more general context of homology cylinders---the
relevant
subgroup is now $\bg$, consisting of those homology cylinders which extend to
homology cylinders over $T_g$.

It turns out that the structure of $\Hg$ presents different problems than $\gg$.
The relative weight filtration extends to a filtration of $\Hg$ but now  the
associated graded group is all of $\D (H)$. On the other hand 
the
residue of the filtration (i.e.\ the intersection of the subgroups of the
filtration) is
non-trivial (it is trivial for $\gg$). There are some approaches to the study of
$\Hg$ which don't seem to have useful analogs for $\gg$. Using the recent work
of Goussarov \cite{Go} and Habiro \cite{H}, there is a notion of finite-type
invariants for homology cylinders. We will show, using results {\em announced}
by Habiro \cite{H}, that this is entirely captured
by the relative weight filtration except for one piece in degree $1$, which is a
$\Z/2$-vector space defined by a natural generalization of the Birman-Craggs
homomorphisms. In a different direction we will exhibit a
close relationship
between
$\Hg$ and a {\em framed string link concordance group} $\sgf$---extending the
natural map
from the
pure braid group into $\gg$ defined by Oda. We will see that $\sgf$ maps  into 
$\Hg$, inducing a {\em bijection} of $\sgf$ with the coset space
$\bhg /\bg$, where $\bhg$ is a natural subset of $\Hg$ containing the
``Torelli'' subgroup---it is interesting to note that there is no analogous
result for $\gg$.

The inclusion $\ggb\sub\gg$ can be illuminated by the introduction of a new
filtration of
$\gg$ which we call the {\em Lagrangian filtration}. The residue of this
filtration is
exactly $\ggb$ and the associated graded group imbeds in a Lie algebra $\D(L)$
constructed from  the
Lagrangian subgroup $L=\ker\{ H_1 (\ssg )\to H_1 (T_g )\}$. Determination of the
image seems as difficult as the analogous problem for the relative weight
filtration. As in the case of the relative weight filtration, the associated
graded group for the Lagrangian filtration of $\Hg$ is isomorphic to $\D(L)$.
However
the residue of the filtration turns out to be {\em larger} than $\bg$.

I would like to thank Stavros Garoufalidis for many useful discussions.

This work was partially supported by NSF grant
DMS-99-71802  and by an Israel-US BSF grant.

\section{Homology cylinders}
\subsection{Preliminaries}

For the usual technical reasons it will be easier to work with the punctured
surface. Let $\Sg$ denote the compact orientable surface of genus $g$  with
 one boundary component. A {\em homology cylinder} over $\Sg$ is a
compact orientable $3$-manifold $M$ equipped with two
imbeddings $i^{-},i^{+}\co\Sg\to\bd M$ so that $i^{+}$ is
orientation-preserving and $i^{-}$ is orientation-reversing and if
we denote $\SS^{\pm}=\im i^{\pm}(\Sg )$, then $\bd M=\SS^+\cup\SS^-$ and
$\SS^{+}\cap\SS^{-}=\bd\SS^{+}=\bd\SS^{-}$. We also require that
$i^{\pm}$ be homology isomorphisms. This notion is introduced in \cite{H}, using
the terminology {\em homology cobordism}. We can multiply two homology
cylinders by identifying $\SS^{-}$ in the first with $\SS^{+}$ in
the second via the appropriate $i^{\pm}$. Thus $\cg$, the set
of orientation-preserving diffeomorphism classes of homology
cylinders over $\Sg$ is a semi-group whose identity is the product
$I\times\Sg$, with $\SS^- =0\times\Sg, \SS^+ =1\times\Sg$, with their collars
stretched half-way along  $I\times\bd\Sg$. $\cg$ is denoted $\C (\ssg )$ in
\cite{H}.

There is a canonical homomorphism $\Gg\to\cg$ that
sends $\phi$ to $(M=I\times\Sg, i^- =0\times\text{id}, i^+ =1\times \phi)$.
Nielsen \cite{N} showed that the natural map $\eta \co\Gg\to \ao (F)$
is an isomorphism, where $F$ is
the free group on $2g$ generators $\{ x_i ,y_i\}$, identified with
the fundamental group of $\Sg$ (with base-point on $\bd\Sg$)---see Figure
\ref{fig1}---and
$\ao  (F)$ is the group of automorphisms of $F$ which fix the
element
$\omega_g =(y_1\cdots y_g)\i (x_1 y_1 x_1\i\cdots x_g y_g x_g\i ) $,
representing
the boundary
of
$\Sg$.

\begin{figure}[ht!]
\psdraw{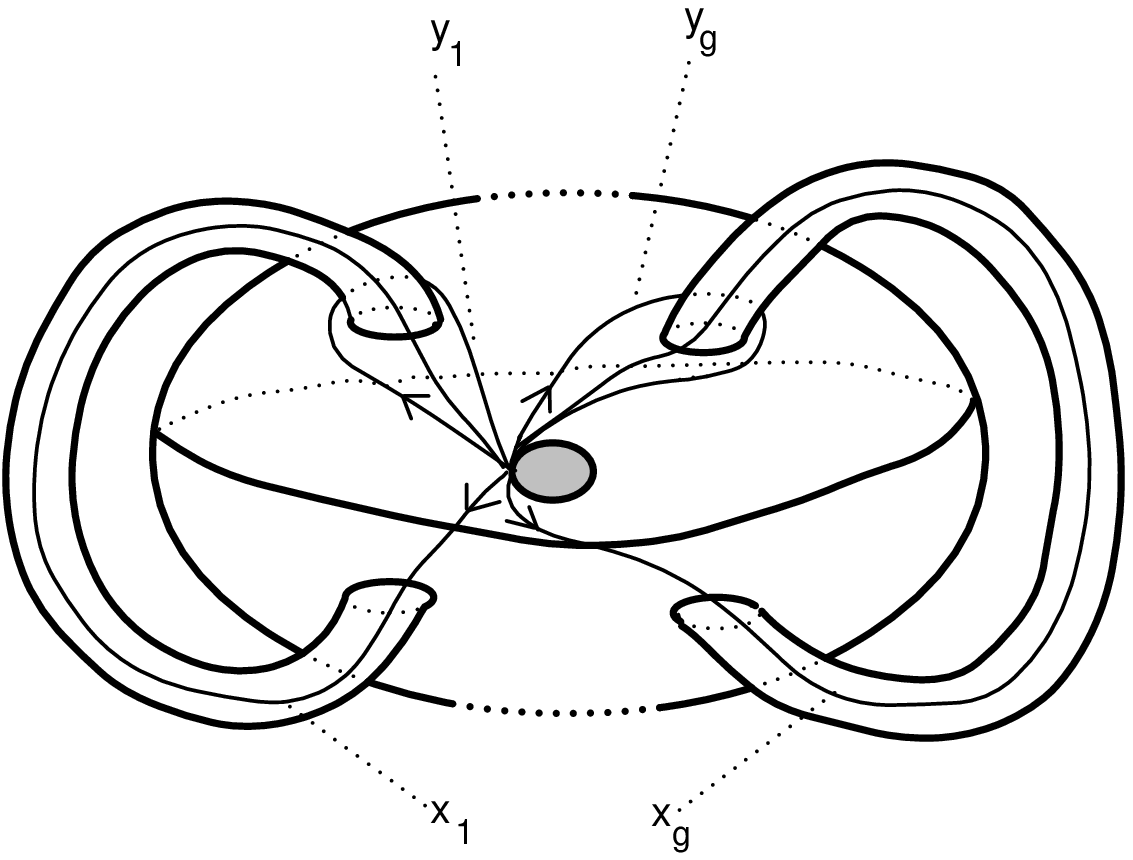}{2.5in}
\caption{Generators of $\pi_1 (\Sg )$}\lbl{fig1}
\end{figure}

We can convert $\cg$ into a group $\Hg$ by considering
{\em homology bordism classes } of homology cylinders. If $M,N$ are homology
cylinders, we can construct a closed manifold $W=M\cup (-N)$, where $\SS^{\pm}$
of
$M$ is attached to $\SS^{\mp}$ of $N$ via their identifications with $\Sg$. A
{\em
homology bordism} between $M$ and $N$ is a manifold $X$ such that $\bd X=W$
and
the inclusions $M\sub X, N\sub X$ are homology equivalences---we say that $M$
and
$N$ are {\em homology bordant}. This is an equivalence relation, since we can
paste
two homology bordisms together to create a third. Furthermore the
multiplication of
homology cylinders preserves homology bordism classes---see Figure \ref{fig3}.
\begin{figure}[ht!]
\psdraw{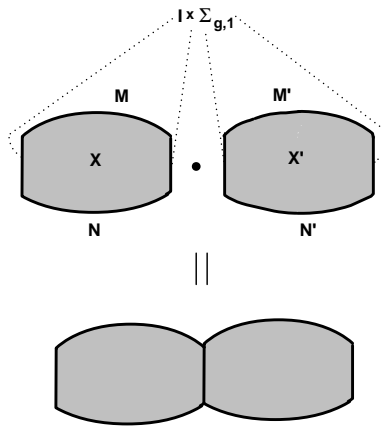}{2in}
\caption{Homology bordism invariance of multiplication of homology
cylinders}\lbl{fig3}
\end{figure}
For any homology cylinder $M$ we can also consider $-M$ as a homology cylinder
with
the roles of $\SS^{\pm}$ reversed. Then $I\times M$ is a homology bordism
between
the product $M\cdot (-M)=\bd I\times M\cup I\times\SS^-$ and $I\times \Sg
=I\times\SS^+$, the identity element of $\cg$. Thus $\Hg$ is a group.

 For any group $G$, let $G_q$ denote the subgroup generated by commutators of
order
$q$. In \cite{GL} the isomorphism $\eta \co\Gg\iso\ao (F)$ is extended to a
sequence of maps 
$$\eta_k \co\Hg\to\ao (F/F_{k+1} )$$
 for $k\ge 1$, which are consistent in the obvious sense. $\ao
(F/F_q
)$ consists of all automorphisms of $F/F_q$ which satisfy the equation 
$$[h(x_1 ),h(y_1 )]\cdots [h(x_g ),h(y_g )]=[x_1 ,y_1 ]\cdots [x_g ,y_g ]\mod
F_{q+1}$$
For example $\ao (F/F_2 )$ is just the symplectic group $\sp (H)$, where $H=H_1
(\Sg )$.

For the reader's convenience we recall the definition of $\eta_k$. Given
$(M,i^{+},i^{-})\in\Hg$ consider
the homomorphisms $i^{\pm}_{*}:F\to\pi_{1} (M)$, where the
base-point is taken in $\bd\SS^{+}=\bd\SS^{-}$. In general, $i^{\pm}_{*}$
are not isomorphisms---however,  since $i^{\pm}$ are
homology isomorphisms, it follows from Stallings \cite{St} that they
  induce isomorphisms
$i^{\pm}_{n}:F/F_{n}\to\pi_{1}(M)/\pi_{1}(M)_{n}$ for any $n$. We then
define $\eta_{k}(M,i^{\pm})=(i^{-}_{k+1})^{-1} \circ i_{k+1}^{+}$. It is
easy to see that $\eta_k(M,i^{\pm}) \in \ao (F/F_{k+1})$.

 One consequence of the existence of these maps is that the map $\Gg\to\Hg$ is
injective. 
Furthermore, the following theorem is proved  in \cite{GL}.

\begin{theorem}\lbl{th.eta}
 Every $\eta_k$ is {\em onto}. 
 \end{theorem}

We now define a filtration of $\Hg$ by setting $\F^w_k (\Hg )=\ker\eta_k$. Then
$\F^w_k  (\Hg )\cap\Gg =\Gg [k]$ is the standard relative weight filtration of
$\Gg$ (this is denoted $\M(k+1)$ in \cite{Mo}). Let $\G_k^w (\Hg )=\F^w_k  (\Hg
)/\F^w_{k +1} (\Hg )$.
There is,
for every $k\ge 1$, a short exact sequence
\begin{equation*}
 1\to\D_k (H)\to\ao (F/F_{k+2})\to\ao (F/F_{k+1} )\to 1 
 \end{equation*}
 where  $\D_k (H)$ is the kernel of the bracketing map $H\otimes\L_{k+1}
(H)\to\L_{k+2}(H)$. $\L_q (H)$ is the degree $q$ part of the free Lie algebra
over
$H$ (see \cite{GL}). 
Theorem~\ref{th.eta} has, as an immediate consequence:
\begin{corollary}
$\eta_{k+1}$ induces an isomorphism $J_k^H \co\G_k^w (\Hg )\iso
\D_k (H)$ for $k\ge 1$. 
\end{corollary}
Note that the induced monomorphism $J_k \co\Gg [k]/\Gg
[k+1]\to\D_k (H)$ is {\em not} generally onto---computing its image is a
fundamental
problem in the study of the mapping class group (see
\cite{J},\cite{Mo},\cite{M}).
We will see that $\F^w_{\infty} (\Hg )=\bigcap_k\F^w_k  (\Hg )$ is non-trivial
whereas it follows from Nielsen's theorem that \newline $\Gg [\infty
]=\{ 1\}$.

\subsection{Filtrations of the Torelli group}
Since $\Gg [1]=\tg$ is the classical Torelli group, we will refer to $\F^w_1 
(\Hg
)=\Tg$ as
the {\em homology Torelli group}. 

It is pointed out in \cite{Mo} that 
$$[\Gg[k],\Gg[l]]\sub\Gg[k+l] $$
Thus  $(\tg )_k\sub\Gg[k]$ . Furthermore the associated graded groups $\G_*^l
(\tg )$
and $\G_*^w (\Gg )$, where $\G_k^l (\tg )=(\tg )_k /(\tg )_{k+1}$ and $\G_k^ w
(\Gg )=\Gg [k]/\Gg [k+1]$, are graded Lie algebras with bracket defined by the
commutator and we have a Lie algebra homomorphism $j\co \G_*^l (\tg )\to\G_*^w
(\Gg
)$ induced by the inclusions. Also $D_* (H)$ is a Lie algebra, as described in
\cite{M},
and the inclusion $J\co\G_*^w (\Gg )\hookrightarrow D_* (H)$ is a Lie algebra
homomorphism. Johnson \cite{J} shows that, after $\otimes\Q$, both $j$ and $J$
are isomorphisms at the degree $1$ level. 
$$ \G_1^l (\tg )\otimes\Q\overset{\iso}\longrightarrow\G_1^ w (\Gg
)\otimes\Q\overset{\iso}\longrightarrow D_1 (H)\otimes\Q$$
Thus $\im j\otimes\Q$ and $\im J\otimes\Q$ are the Lie
subalgebras of $D_* (H)\otimes\Q$ generated by elements of degree $1$.

The exact same considerations apply to the lower central series filtration of
$\Hg$ and
we obtain a Lie algebra homomorphism
$$j^H \co\G_*^l (\Tg )\to\G_*^w (\Hg )\iso D_* (H) $$
The image of $j^H\otimes\Q$ is again the subalgebra generated by degree $1$
elements,
which is the same as $\im J\otimes\Q$, and so we conclude that $j^H$ is not
onto. 
\subsection{Finite-type invariants of $\Hg$}
  There have been several proposals for a theory of finite-type invariants for
general $3$-manifolds, extending Ohtsuki's theory \cite{Oh} for homology
$3$-spheres---perhaps the first given in Cochran-Melvin \cite{CM}. We will use,
however,
the particular version proposed independently by M. Goussarov
\cite{Go},\cite{GGP} and K. Habiro \cite{H}. We will show that $\eta_k$ and a
homomorphism $\b \co\Tg\to V$ defined by Birman-Craggs homomorphisms \cite{BC},
where $V$ is a vector space over $Z/2$, make up the
{\em universal multiplicative homology bordism invariant of type $k$} on the
class of homology cylinders for any $k\ge 1$. It is to be emphasized that {\em
we make use of
results announced in \cite{H} but, for which no proofs have yet appeared.}

We give a brief summary of the Goussarov-Habiro theory, and refer the reader to
\cite{GGP},\cite{H},\cite{Ha} for  details. 
Let $G$ be a unitrivalent graph whose trivalent vertices are equipped with a
cyclic
ordering of its incident edges and whose univalent vertices are decorated with
an element
of an abelian group $H$. We also insist that each component of $G$ have at least
one trivalent vertex. We refer to such a graph as an {\em $H$-graph} and
define the
{\em degree} of $G$ to be the number of trivalent vertices. If $M$ is a
$3$-manifold
and $H=H_1 (M)$, then a {\em clasper} (using the terminology of \cite{H}---{\em
clover} in the terminology of \cite{GGP}) in $M$
associated to $G$ is a framed link
$C$ in
$M$ obtained in the following way. Associate to each trivalent vertex of $G$ a
copy of the Borromean rings, in disjoint balls
of $M$, each component associated to an end of an edge incident to that vertex
(even if two of those ends are just opposite ends of the same edge) and given
the $0$-framing. Associate
to each univalent vertex a framed knot in $M$, disjoint from the balls, and
representing the element of $H$ labeling that vertex. These are the {\em leaves}
of the clasper. Finally for each edge of $G$ introduce a simple clasp between
the knots
associated to the two ends of that edge.  The construction of a clasper $C$ from
an $H$-graph $G$ involves a
choice of framed  imbedding into $M$ of the graph obtained from $G$ by attaching
circles to the univalent vertices. See
Figure \ref{fig.clasper} for a typical example. 
\begin{figure}[ht!]
\psdraw{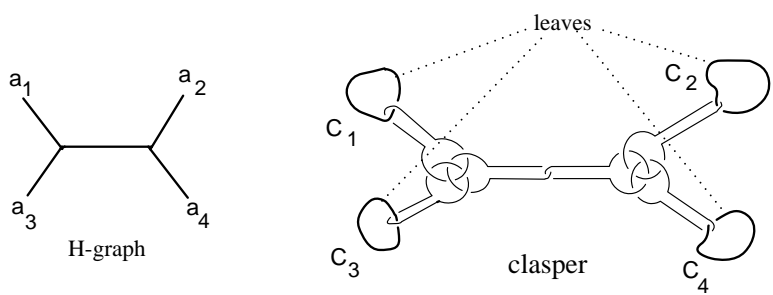}{4in}
\caption{An H-graph and an associated clasper, where $C_i$ represents
$a_i$}\lbl{fig.clasper}
\end{figure}
If $G$ has degree $n$ we say $C$ is an
{\em $n$-clasper}. Surgery on a clasper produces a new manifold $M'=M_C$ with
isomorphic
homology and torsion pairing. We will refer to surgery on an $n$-clasper as an
{\em $n$-surgery}.
Thus if $M$ is a homology cylinder of genus $g$ so
is $M'$. It is easy to see that the associated automorphisms of $H_1 (\Sg )$
associated to $M$ and $M'$ are the same. Matveev \cite{Ma} has proved  the
converse: If $M, M'$ have isomorphic homology and torsion pairing then there is
a clasper $C$ in $M$ such that 
$M'=M_C$. 
We define the relation of {\em $A_k$-equivalence} (in the terminology of
\cite{H}) to be generated by the following elementary move: $M\sim_k M$ if $M'$
is diffeomorphic to $M_C$, for some connected $k$-clasper $C$. According to
\cite{H}, the set of $A_k$-equivalence classes in $\cg$ is a group
(finitely-generated and nilpotent) under the stacking multiplication. Let $\F_k
(\cg )$ denote the subsemigroup of $\cg$ consisting of homology cylinders
$A_k$-equivalent to the trivial one $I\times\Sg$ and let $\G_k (\cg )=\F_k (\cg
)/A_{k+1}$

Habiro defines $\A_k (H)$, for $k>1$, to be the abelian group generated by
connected $H$-graphs of
degree $k$ with an
extra structure of a total ordering of the univalent vertices, subject to the AS
and IHX relations, multilinearity of the labels and an STU-like relation(see
\cite{H}). $\A_1 (H)$ is defined explicitly to be $\LL^3 H\oplus V$, where
$V=\LL^2 H_2\oplus H_2\oplus\Z/2$, where $H_2=H\otimes\Z/2$. Habiro constructs
an
{\em epimorphism} $\A_k (H)\to\G_k (\cg )$, which can be defined for $k>1$ by
clasper surgery 
on $I\times\Sg$, in 
which the ordering of the univalent vertices is used to stack the imbeddings
(with horizontal framings) of the leaves 
representing
the elements of $H=H_1 (\Sg )$ decorating those vertices. Habiro conjectures
that this epimorphism is an isomorphism and claims
that it is so for $k=1$ and, for $k>1$,  induces an isomorphism $\A_k
(H)\otimes\Q\iso\G_k (\cg
)\otimes\Q$. See \cite{Ha} for an alternative construction for homology
handlebodies.

Let $\A^{\circ}_k (H)$ denote, for $k>1$, the subgroup of $\A_k (H)$ generated
by connected
$H$-graphs with
non-zero first Betti number, i.e.\ with a non-trivial cycle. The quotient $A^t_k
(H)=\A_k (H)/\A^{\circ}_k (H)$ is
generated by {\em $H$-trees}. Note that the ordering of the univalent vertices
in $\A^t_k (H)$ does not matter, because of the STU-relation, and so this
structure can be ignored. Set $\A^t_1 (H)=\A_1 (H)$.

We recall a fact known to the experts (see also \cite{Ha}).

\begin{proposition}\lbl{prop.tree}
$\A^t_k (H)\iso\D_k (H)$ for $k>1$. $\A^t_1 (H)\twoheadrightarrow\D_1
(H)\iso\LL^3 H$.
\end{proposition}
\begin{proof}[Outline of proof.]
If $k=1$, this was already pointed out by Johnson \cite{J}. Now suppose $k>1$.
Let $T$ be an $H$-tree of degree $k$. If we choose any univalent vertex
to be the root of $T$, then we can associate to this rooted $H$-tree an element
in $H\otimes L_{k+1}(H)$. For example, if
$T$ is the $H$-graph in Figure \ref{fig.clasper}, with the counter-clockwise
orientation at each vertex and root at $a_1$,
then
the associated element would be $a_1\otimes [\a_3 ,[\a_4 ,\a_3 ]]$, where $\a_i$
is a lift of $a_i$ into $\pi_1 (\Sg )$.
Now
we associate to $T$  the sum of the associated elements of $H\otimes L_{k+1}(H)$
 as the roots range over all the univalent vertices  of $T$. It is not
hard to see
that this element actually lies in $\D_k (H)$ and that, in fact,  $\D_k (H)$ is
generated by the elements associated to all
possible
$H$-trees of degree $k$.
\end{proof}

To bring homology bordism into the picture we need the following theorem.

\begin{theorem}\lbl{th.homcob}
Let $M$ be a $3$-manifold and $G$ a connected $H$-graph, where $H=H_1 (M)$, with
at least
one non-trivial cycle. If $C$ is any clasper associated to $G$, then $M_C$ is
homology bordant to $M$.
\end{theorem}
 \begin{remark}\lbl{rem.link}
 As an immediate consequence of this theorem we conclude that if $L\sub M$ is a
link in a $3$-manifold and $C$ is such a clasper in $M-L$ then  the links
$(M,L)$ and $(M_C ,L)$ are {\em homology concordant}, i.e.\ there is a homology
bordism $V$ between $M$ and $M_C$ and a proper imbedding $I\times L\sub V$ such
that $0\times L=L\sub M$ and $1\times L=L\sub M_C$.
 Compare this to Theorem 2.9 of \cite{Ha} which shows that, in the case of $M$ a
homology ball, the Milnor $\bar\mu$-invariants of $L\sub M$ and $L\sub M_C$
coincide.
\end{remark}
\begin{remark}
A more delicate argument will strengthen Remark \ref{rem.link} as follows. If 
the clasper $C$ is
{\em strict} in the sense of \cite{H}, i.e.\ the leaves of $C$ bound disks which
are disjoint from each other and the rest of $C$ except for a single
intersection point with the companion component of $C$ (but will generally
intersect $L$), then $(M_C ,L)$ is {\em concordant} to $(M,L)$, i.e.\ there is an
imbedding $I\times L\sub I\times M$ so that $(0\times M,0\times L)$ is
diffeomorphic to $(M,L)$ and $(1\times M,1\times L)$ is diffeomorphic to $(M_C
,L)$.
\end{remark}
To prove Theorem \ref{th.homcob} we need the following Lemma.
\begin{lemma}\lbl{lem.homcob} Let $M$ be a $3$-manifold and $L\sub M$ a framed
link. Suppose
that $L=L'\cup L''$, where $L'$ and $L''$ have the same number of components and 
\begin{itemize}
\item[\rm(a)] $L'$ is a trivial link bounding disjoint disks $D$,
\item[\rm(b)] $L'$ is $0$-framed,
\item[\rm(c)] the matrix of intersection numbers of the components of $D$ with
those
of
$L''$ is non-singular.
\end{itemize}
Then $M_L$ is homology bordant to $M$.
\end{lemma}
\begin{proof} We construct a manifold $V$ from $I\times M$ by adjoining
handles to $0\times M$ along the framed link $L''$ and by removing tubular
neighborhoods of properly imbedded disjoint disks $D'$ in $I\times M$ obtained
by pushing $0\times\inte D$ into the interior of $I\times M$. Then $\bd V=M_L
-(1\times
M)$, so it only remains
to observe that the pair $(V, 1\times M)$ is acyclic.

Let $W=(I\times M)-D'$. It is easy to see that the only non-zero homology group
of
$(W,1\times M)$ is $H_1 (W,1\times M)$, which is freely generated by the
meridians of $D'$. By considering the triple $(V,W,1\times M)$ we find there is
an exact sequence:
$$0\to H_2 (V,1\times M)\to H_2 (V,W)\to H_1 (W,1\times M)\to H_1 (V,1\times
M)\to 0 $$
Now $H_2 (V,W)$ is the only non-zero homology group of $(V,W)$ and it is freely
generated by the disks adjoined along $L''$. Since the homomorphism $H_2
(V,W)$ $\to H_1 (W,1\times M)$ is represented by the matrix of intersection numbers
of the components of $D'$ with those of $I\times L''$, it follows that
$(V,1\times M)$ is acyclic if and only if (c) is satisfied.
\end{proof}
\begin{proof}[Proof of Theorem~\ref{th.homcob}]
 Recall that for each edge of $G$ there are two components of $C$, one at each
end of the edge (see Figure~\ref{fig.clasper}). We will call them {\em
companion} components. We will construct $C'$ and $C''$ by assigning, for each
edge of $G$, one of the associated companions to $C'$ and the other to $C''$.
This choice can be represented by an orientation of the edge pointing from the
end associated to the companion in $C'$ toward the end associated to $C''$. Our
aim will be to make these choices satisfy:
\begin{itemize}
\item[(i)] An edge with a univalent vertex ( a {\em leaf} of $G$) is oriented
toward the leaf, i.e.\ {\em outward} (thus the leaves of $C$ will all belong to
$C''$),
\item[(ii)] No trivalent vertex is a source, i.e.\ not all the incident edges are
oriented away from the vertex.
\end{itemize}
We will see that these conditions can be satisfied if and only if $G$ is not a
tree. But for now note that if these conditions are satisfied then the
decomposition of $C$ will satisfy conditions (a)-(c) of Lemma~\ref{lem.homcob}. 
Since the three components of $C$ associated to any trivalent vertex are a
Borromean rings, any two of the components bound disjoint disks. Thus we can
choose disjoint disks $D$ bounded by each component of $C'$, the disks from
components associated to different vertices will be disjoint. The only problem
would be if these components were associated to the same edge, but this is ruled
out. Now each disk from a component of $C'$ will intersect the companion
component of $C$ once---the only other intersections will be with one of the
components of $C''$ which is associated to the same trivalent vertex, but the
intersection number will be $0$. Thus the intersection matrix of the components
of $D$ with those of $C''$ will be the identity matrix. 
See Figure~\ref{fig.dcom} for an example. 
\begin{figure}[ht!]
\psdraw{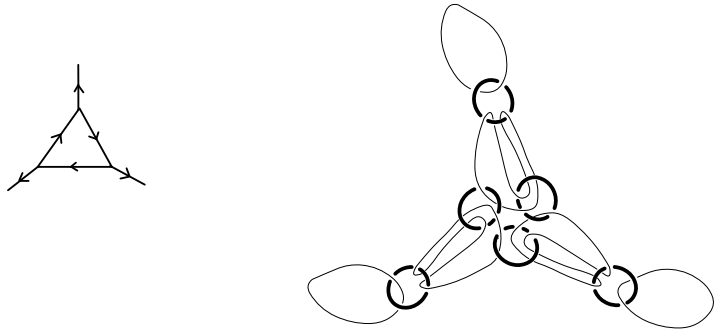}{3.7in}\vglue -15pt
\caption{An oriented graph on the left and a corresponding clasper $C$ on the
right. The bold components define $C'$.}\lbl{fig.dcom}
\end{figure}

\begin{lemma}\lbl{lem.ori}
 Suppose that $T$ is a unitrivalent tree, with an orientation prescribed for
each leaf edge, not all outward. Then we can extend this to an orientation of
all the
edges of $T$ which satisfies (ii).
\end{lemma}
\begin{proof} A similar fact is proved in \cite{Ha}. Choose one of the leaves
$e$ of $T$ which is oriented inward. Now
orient every edge of $T$ which is not a leaf so that it points {\em away} from
$e$, i.e.\ if we travel along any non-singular edge path which begins at $e$ and
ends at a non-leaf, then the orientations of all the edges in the path point in
the direction of travel. Then
it is clear that any trivalent vertex will have at least one of its
incident edges oriented toward that vertex. Thus (ii) will be satisfied whatever
the orientations of the other leaves. See Figure \ref{fig.ex} for an
example. \end{proof}

\begin{figure}[ht!]
\psdraw{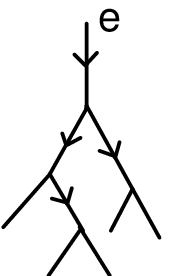}{0.55in}
\nocolon\caption{}\lbl{fig.ex}
\end{figure}

We can now complete the proof of Theorem~\ref{th.homcob}. Since $G$ has a cycle
we can make one or more cuts in edges of $G$ to create a tree. Each edge of
$G$ which is cut will create two new leaves in $T$. We now choose arbitrary
orientations of
each cut edge of $G$, which will induce orientations of the new leaf edges of
$T$.
Note that one of each pair of new leaves will be oriented inwards. Thus the
outward orientations of the leaves of $G$ together with these orientations of
the new leaves of $T$ provides orientations of all the leaves of $T$ which
satisfies the hypothesis of Lemma~\ref{lem.ori}. Applying this lemma gives an
orientation of $T$ satisfying (ii). But now we can glue the cut edges back
together and we get an orientation of $G$ satisfying (i) and (ii), thus proving
the Theorem. See Figure~\ref{fig.cut} for an example.\end{proof}
\begin{figure}[ht!]
\psdraw{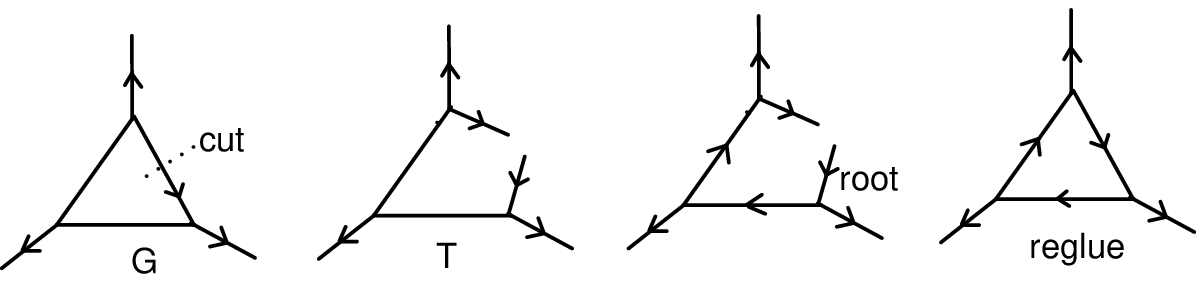}{4.2in}
\nocolon\caption{}\lbl{fig.cut}
\end{figure}

Let's define a filtration of $\Hg$ by $\F_k^Y (\Hg )=\im (\F_k (\cg )\to\Hg )$,
where we consider the restriction of the quotient map $\cg\to\Hg$. By \cite{Ma},
$\F_1^Y (\Hg )=\F_1^w (\Hg )=\Tg$.

\begin{theorem}\lbl{th.univ}
$\F_k^Y (\Hg )\sub\F_k^w (\Hg )$, thus inducing a map $\G_k^Y (\Hg )\to\G_k^w
(\Hg )$ and a commutative diagram:
\begin{equation}\lbl{eq.univ}
\begin{diagram}
\node{\A_k (H)}\arrow{s}\arrow{e}\node{\G_k (\cg )}\arrow{s}\\
\node{\A_k^t (H)}\arrow{e,t}{\theta_k}\node{\G_k^Y (\Hg )}\arrow{e}\node{\G_k^w
(\Hg
)}\arrow{e,tb}{J^H_k}{\iso}\node{\D_k (H)}
\end{diagram}
\end{equation}
$\theta_k$ is an isomorphism for all $k$ and the composition
$\A_k^t (H)\to\G_k^Y(\Hg )\to\G_k^w (\Hg )\to\D_k (H) $
is the epimorphism (isomorphism if $k>1$) of Proposition~\ref{prop.tree}.

\end{theorem}
\begin{proof}

First note that it follows from Theorem~\ref{th.homcob} that the composition
$\A_k (H)$ $\to\G_k (\cg )\to\G_k^Y (\Hg )$ factors through $\A^t_k (H)$. 
This yields the commutative diagram~\eqref{eq.univ}, assuming for the moment
that $\F_k^Y (\Hg )\sub\F_k^w (\Hg )$.

To see that $\F_k^Y (\Hg )\sub\F_k^w (\Hg )$ for $k>1$, and, at the same time,
identify the
composition of the maps in the bottom line of the diagram, we
first consider a $k$-clasper $C$ in $I\times\Sg$ associated to a connected
$H$-tree $T$. If we choose
a root of $T$ and the element of $H\otimes
L_{k+1}(H)$ associated to this rooted $H$-tree is $\a\otimes \g$, then, by a
sequence of Kirby moves, we can convert $C$
into a
$2$-component link $(l,l')$, where $l$ is the leaf of $C$ corresponding to the
root, and so represents $\a\in\pi / \pi_2$,
and
$l'$ represents $\g\in\pi_{k+1} \mod\pi_{k+2}$ . For example, Figure
\ref{fig.comm} explains this for a rooted $H$-tree
of
degree $1$.

\begin{figure}[ht!]
\psdraw{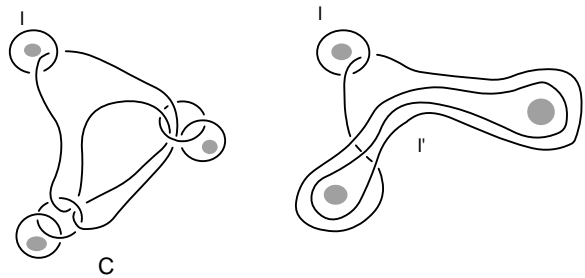}{3in}
\caption{A clasper $C$ on the left, and the Kirby-equivalent link on the
right}\lbl{fig.comm}
\end{figure}

Now suppose $M=(I\times\Sg )_C$ and $\xi\in\pi_1 (\SS_+ )$. To compute $\eta_k 
(M)\cdot\xi$ we push a curve $K$
representing $\xi$ down $I\times\Sg$ and observe how it changes every time we
cross $C$. When it crosses the leaf $l$, the
effect may be computed by using the Kirby-equivalent link $l\cup l'$ instead of
$C$. But then it is not hard to see that the
effect is to add a copy of $l'$ to $K$. Thus the contribution of all the
crossings of $K$ with the leaf $l$ is $<\xi ,\a >\g$, where
$<,>$ denotes the homological intersection number in $\Sg$. Under the canonical
isomorphism $\hom (H,L_{k+1}(H))\iso H\otimes
L_{k+1}(H)$ this corresponds to $\a\otimes\g$. The total change in $\xi$ is then
the sum of these contributions over all the
univalent vertices of $T$, which is exactly the element of $\D_k (H)$
corresponding to $T$. This shows that $\F_k^Y (\Hg )\sub\F_k^w (\Hg )$, at least
for the elements represented by trees in  $\F_k^Y (\Hg )$, and that the
composition $\A_k (H)\to\G_k (\cg )\to\G_k^Y (\Hg )\to\G^w_k (\Hg )\to\D_k (H)$
associates to any $H$-tree the corresponding element of $\D_k (H)$ under the
isomorphism of Proposition~\ref{prop.tree}. The proof of the commutativity of
diagram~\eqref{eq.univ} is completed by noting that if $G$ is any $H$-graph with
a cycle
(and an ordering of its univalent vertices), then we can cut open some edges to
create a tree $T$, where the new vertices are labeled by $0\in H$. Now a clasper
representing $G$ is also a clasper representing $T$ and the above argument
applies. Since some of the labels are $0$, the image in $\D_k (H)$ is $0$.

To prove that $\theta_k$ is an isomorphism first note that the map $\A_k^t
(H)\to\G_k^Y
(\Hg )$ is onto, since $\A_k (H)\to\G_k (\cg )$ is onto. So if $k>1$ it must
also be
one-one, since the composition to $\D_k (H)$ is an isomorphism. Thus all the
maps
in the bottom row must be isomorphisms.
\end{proof}

To deal more completely with $k=1$ we must consider some additional
$\Z/2$-invariants on
$\Hg$---or
rather on $\Tg =\F_1^Y (\Hg )=\F_1^w (\Hg )$. Recall the Birman-Craggs
homomorphisms on the classic Torelli group $\tg$ from \cite{BC}. To define these
first
choose an imbedding $i\co\ssg\sub S^3$. Then if $h\in\tg$ we can cut $S^3$ along
$i(\ssg )$ and
reglue, using $h$, to obtain a homology sphere $\SS_h$. We then define $\b_i
(h)$ to be
the Rochlin invariant of $\SS_h$. It is clear that this makes perfectly good
sense if we
substitute a homology cylinder in $\Tg$ for $h$ and the result depends only on
its
homology bordism class. Thus we obtain homomorphisms $\b_i \co\Tg\to\Z/2$. 

Johnson shows that one can choose a number of $\{\b_i\}$ exactly equal to the
dimension
of the vector space $V$, thereby defining a homomorphism $\b \co\tg\to V$, 
which
determine all the Birman-Craggs homomorphisms, and he then shows that the
combined
homomorphism $H_1 (\tg)\to\LL^3 H\oplus V$ is an isomorphism.

Now consider the corresponding $\b\co\Tg\to V$. We claim that $\b (\F_2^Y (\Hg
))=0$ and
so we get an induced $\b \co\G^Y_1 (\Hg )\to V$. As is pointed out in \cite{H}
and
\cite{GGP},
surgery on a clasper of degree $n$ is the same as cutting and pasting along some
imbedded
surface using an element of the $n$-th lower central series term $(\tg )_n$ of
the
Torelli group. Thus a Birman-Craggs homomorphism on a homology cylinder in
$\F_n^Y
(\Hg )$ takes the same value as some Birman-Craggs homomorphism on some element
of
$(\tg )_n$, which is always zero if $n\ge 2$. 

Let us now consider the combined homomorphism $\G^Y_1 (\Hg )\to\G^w_1 (\Hg
)\oplus
V\iso\LL^3 H\oplus V$, which must be an epimorphism since it is when restricted
to
$\tg$. Since we also have Habiro's epimorphism $\LL^3 H\oplus V\iso\A_1^t
(H)\to\G_1^Y (\Hg )$, it follows that both of these epimorphisms are actually
isomorphisms.

To summarize, we have:
\begin{corollary}\lbl{cor.fti} 
\begin{enumerate}
\item $\G^Y_1 (\Hg )\iso\G^w_1 (\Hg
)\oplus
V\iso\LL^3 H\oplus V$,
\item $\G_k^Y (\Hg )\iso\G^w_k (\Hg )$ for $k>1$,
\item 
$\F_k^Y (\Hg )=\F_k^w (\Hg )\cap\ker\b$ for $k>1$.
\end{enumerate}
\end{corollary}
\begin{proof} (1) and (2) are proved above. To prove (3) first consider $k=2$.
Suppose that $M\in\F_2^w (\Hg )\sub\F_1^w (\Hg )=\F_1^Y (\Hg )$ and $\b (M)=0$.
Then $M\to 0\in\G_1^w (\Hg )\oplus V\iso\G_1^Y (\Hg )$ and so $M\in\F_2^Y (\Hg
)$. For $k\ge 3$ we proceed by induction, using the following commutative
diagram.
$$
\begin{diagram}
\node{}\node{0}\arrow{s}\node{0}\arrow{s}\\
\node{0}\arrow{e}\node{\F_k^Y (\Hg )}\arrow{s}\arrow{e}\node{\F_k^w (\Hg
)}\arrow{s}\arrow{e,t}{\b}\node{V}\arrow{s,=}\\
\node{0}\arrow{e}\node{\F_{k-1}^Y (\Hg )}\arrow{e}\arrow{s}\node{\F_{k-1}^w (\Hg
)}\arrow{s}\arrow{e,t}{\b}\node{V} \\
\node{}\node{\G_{k-1}^Y (\Hg )}\arrow{s}\arrow{e,t}{\iso}\node{\G_{k-1}^w (\Hg
)}\arrow{s}\\
\node{}\node{0}\node{0}
\end{diagram}
$$
The columns are exact and, by induction, the middle row is exact. It then
follows that the top row is exact.
\end{proof}

Suppose that $\phi$ is a multiplicative invariant of homology cylinders of genus
$g$, i.e.\ a homomorphism $\phi \co\cg\to G$ for some group $G$. We will say that
$\phi$ is of {\em finite type} if $\phi (\F_{k+1}(\cg ))=1$ for some $k$, and is
of type $k$ for the minimum such value. This is actually a {\em weaker}
definition than the usual notion of finite type. 

For example $\eta_k$ is an invariant of type $k$, as is $\b$, although it is
only defined on $\Tg$. The following corollary
asserts that $\eta_k$ and $\b$ are the {\em universal} homology bordism
invariant of type $k$.

\begin{corollary} Suppose $\phi$ is a multiplicative homology bordism invariant
of
type $k$. If $M, N$ are two homology cylinders such that $\eta_k (M)=\eta_k (N)$
and $\b (M\cdot N\i )=0$, then $\phi (M)=\phi (N)$. 
\end{corollary}
\begin{proof} By Corollary \ref{cor.fti} $M\cdot N\i\in\F_{k+1}^Y (\Hg )$ and so
$\phi (M\cdot N\i )=1$.
\end{proof}

\section{String links and homology cylinders}
\subsection{String links}\lbl{sec.str}
Let $\P_g$ denote the group of pure braids with $g$ strands. Let $\S_g$ denote
the
group of concordance classes of string links with $g$ strands in a homology
$3$-ball
whose boundary is identified with $\bd (I\times D^2 )$. Two string links $S_1
,S_2$
in homology $3$-balls $B_1 ,B_2$, respectively, are {\em concordant} if there is
a
homology $4$-ball $B$ whose boundary is identified with $B_1\cup B_2$  (with
their
boundaries identified), and a concordance between $S_1$ and $S_2$ imbedded in
$B$.
There is an obvious homomorphism $\P_g\to\S_g$. Recall the theorem of Artin
which states
that
the map $\pg\to\aut (F')$---where $F'$ is the free group with basis $\{ y_1
,\dots
,y_g\}$ identified with $\pi_1 (\dg )$, and $\dg$ is the $2$-disk with $g$
holes---is
injective with image $\ai (F')$. $\ai (F')$ is the subgroup of $\aut (F')$
consisting of
all automorphisms $h$ such that $h(y_i )=\l_i\i y_i\l_i$, for some choice of
$\l_i$
and $h(y_1\cdots y_g )=y_1\cdots y_g$. Note that $\l_i$ is uniquely determined
if we
specify that the exponent sum of $y_i$ in $\l_i$ is $0$. 

The Milnor $\bar\mu$-invariants can be formulated as a sequence of homomorphisms
$\bar\mu_k \co\sg\to\ai (F'/F'_{k+2} ),\ k\ge 1$, where $\ai (F'/F'_q )$
consists
of all
automorphisms $h$ satisfying $h(y_i )=\l_i\i y_i\l_i$ for some $\l_i\in
F'/F'_{q-1}$ and satisfying the equation 
\begin{equation}\lbl{eq.mu}
h(y_1\cdots y_g )=y_1\cdots y_g
\end{equation}
 (see \cite{HL}).  Note now that $\l_i$ is uniquely determined by $h$ up to left
multiplication by a power of $y_i$. One consequence of the existence of the
$\bar\mu_k$ is that the map $\pg\to\sg$ is injective. 
It is known that $\bar\mu_k$ is onto (see e.g.\ \cite{HL}). If we define $\sg
[k]=\ker\bar\mu_k$ and $\sg [0]=\sg$, then we have an isomorphism $\sg [k]/\sg
[k+1]\iso \D_k
(H')$, for $k\ge 0$, where $H'=H_1 (\dg )$. This isomorphism is defined by
$\s\to
\sum_i y_i\otimes\l_i$, for $\s\in\sg [k]$, where $\l_i\in F_{k+1}/F_{k+2}$ are
determined by $\mu_{k+1}$. Note that $\sg [\infty ]=\bigcap_k\sg [k]$ is
non-trivial
since, for example, it contains the knot concordance group.

If $\sgo$ denotes the standard string link concordance group, i.e.\ string links
in
$I\times D^2$ and concordances in $I\times (I\times D^2 )$, then we have natural
maps $\pg\sub\sgo\to\sg$. 

\begin{question} Is $\sgo\to\sg$ injective? Since the realization theorem of
\cite{HL}
produces standard string links the induced map $\sgo /\sgo [\infty ]\to\sg /\sg
[\infty
]$ is an isomorphism.
\end{question}

If $\th$ denotes the group of homology bordism classes of closed homology
$3$-spheres, then
we have an obvious injection $\th\sub\sg$, defined by connected sum with the
trivial string link in $I\times D^2$, whose image is a central subgroup of
$\sg$, lying in $\sg
[\infty
]$, and a retraction $\sg\to\th$. 

\begin{question} Is the combined map $\th\times\sgo\to\sg$ an isomorphism?
\end{question}

We will need to consider the {\em framed} versions. Let $\pgf$ and $\sgf$ denote
the
groups of framed pure braids on $g$ strands and concordance classes of framed
string
links with $g$ strands. Note that $\pgf\iso\pg\times\Z^g$ and
$\sgf\iso\sg\times\Z^g$.
We have an isomorphism $\mu \co\pgf\to\A_1 (F')$ and homomorphisms $\mu_k
\co\sgf\to\A_1 (F'/F'_{k+1} )$ for $k\ge 1$. $\A_1 (F')$ (resp.\ $\A_1 (F'/F'_q
)$)
consists of
sequences $\l =(\l_1 
,\dots
,\l_g )$ of elements $\l_i\in F'$ (resp.\ $F'/F'_{q}$) which satisfy the
following

\begin{enumerate}
\item[(i)] The map $y_i\to\l_i\i y_i\l_i$ defines an {\em automorphism
}$\phi_{\l}$ of $F'$ (resp.\ $F'/F'_{q+1}$)
\item[(ii)] $\phi_{\l} (y_i\cdots y_g )=y_1\cdots y_g$.
\end{enumerate}

The multiplication in $\A_1 (F'), \A_1 (F'/F'_q )$ is defined by $(\l\mu )_i
=\l_i\phi_{\l}(\mu_i )$. Then $\l\to\phi_{\l}$ defines  epimorphisms
with 
kernel
$\Z^g$. Note that we have the commutative diagram
$$
\begin{diagram}
\node{\sgf}\arrow{e}\arrow{s,l}{\mu_k}\node{\sg}\arrow{s,r}{\bar\mu_k}\\
\node{\A_1 (F'/F'_{k+1} )}\arrow{e}\node{\ai (F'/F'_{k+2 })}
\end{diagram}
$$
We define the filtration $\sgf [k]=\ker\mu_k$.

\subsection{Relating string links and homology cylinders}\lbl{sec.sl}

Recall the imbedding $\Phi \co\pgf\to\Gg$, defined by \cite{O} and studied in
\cite[Section 2.2]{L}. There
is a commutative diagram
$$
\begin{diagram}
\node{\pgf}\arrow{e,t}{\Phi}\arrow{s,l}{ \mu}\node{\Gg}\arrow{s,r}{\eta}\\
\node{\A_1 (F')}\arrow{e,t}{\phi}\node{\ao (F)}
\end{diagram}
$$
where $\phi$ is defined by 
\begin{equation}\lbl{eq.phi}
\begin{split}
\phi ((\l_i ))\cdot y_i &=\l_i\i y_i\l_i \\
\phi ((\l_i ))\cdot x_i &= x_i\l_i
\end{split}
\end{equation}
Since $\mu$ and $\eta$ are isomorphisms and $\phi$ is injective, it follows that
$\Phi$ is injective.  

We now extend $\Phi$ to an imbedding $\hp \co\sgf\to\Hg$. Choose an imbedding
$\dg\sub\Sg$ where the meridians of the holes in $\dg$, which correspond to the
generators $y_i\in F'$, are mapped to the meridians of the handles of $\Sg$
which
correspond to the generators $y_i\in F$---see Figure \ref{fig1}. Suppose $S$ is
a framed string link
with $g$
strands in a homology $3$-ball $B$ and $C$ is the complement of an open tubular
neighborhood of $S$ in $B$. Then there is a canonical identification of $\bd C$
with
$\bd (I\times\dg )$. Now we create a homology cylinder by taking $I\times\Sg$,
with $\SS^- =0\times\Sg,\ \SS^+ =1\times\Sg$,
removing $I\times\dg$  and replacing it with $C$, using the canonical
identification of the boundaries. We take this homology cylinder to represent
$\hp (S)$.

It is clear that $\hp$ is a well-defined homomorphism and extends $\Phi$ (see
the
definition of $\Phi$ in \cite{L}). By examining Figure \ref{fig1} we can see
that the following diagram commutes:
$$
\begin{diagram}
\node{\sgf}\arrow{e,t}{\hp}\arrow{s,l}{ \mu_k}\node{\Hg}\arrow{s,r}{\eta_k}\\
\node{\A_1 (F'/F'_{k+1})}\arrow{e,t}{\hat\phi}\node{\ao (F/F_{k+1})}
\end{diagram}
$$
where $\hat\phi$ is the injection defined by 
\begin{equation}\lbl{eq.hphi}
\begin{split}
\hat\phi ((\l_i ))\cdot y_i &\con\l_i\i  y_i\l_i \mod F_k\\
\hat\phi ((\l_i ))\cdot x_i &\con x_i\l_i \mod F_k
\end{split}
\end{equation}
As a consequence we see that $\hat\Phi$ preserves the weight
filtrations and the induced map
$$\sgf [k]/\sgf [k+1]\longrightarrow\F^w_k  (\Hg )/\F^w_{k+1}(\Hg )=\G^w_k
(\Hg ) $$
corresponds to the {\em monomorphism} $\D_k (H')\to\D_k (H)$ induced by
inclusion
$H'\sub H$.

\begin{theorem}\lbl{th.Phi}
 $\hp$ is injective.
\end{theorem}

\begin{proof}
Let $\bhg$ denote the subset (it is not a subgroup) of $\Hg$ consisting of all
$M$ such
that $p\eta_1 (M)i$ is an isomorphism, where $i\co H'\to H$ is the inclusion and
$p\co H\to
H'$ is the projection with $p(x_i )=1,\ p(y_i )=y_i$. Note that $\hp (\sgf
)\sub\bhg$. We will define a
``retraction'' $\rho
\co\bhg\to\sgf$, which is not a homomorphism, but will satisfy $\rho\circ\hp
=\text{id}$
and so prove the theorem. 

We identify $I\times\dg$ with the complement of the trivial framed string link.
Now
consider an imbedding of $\dg$ into $\frac12\times\dg$ defined by removing a
thin
collar of the entire boundary of $\frac12\times\dg$. Thus we obtain an imbedding
of $\dg$ into the
interior
of $I\times\dg$. Now let $\hg$ denote the solid handlebody of genus $g$ whose
boundary
is $\ssg\supseteq\Sg$. We will make this identification so that
the
$y_i$ represent a basis for $\pi_1 (\hg )$, which we identify with $F'$, and the
$x_i$
correspond to the meridians of the handles of $\hg$. We can imbed $\hg$ into
$I\times\dg$ as a thickening of the imbedded $\dg$ and so that
$\dg\sub\Sg\sub\hg$ is
the
imbedding which we used to define $\hp$.

Now suppose $M$ represents an element of $\Hg$. Then we can cut open
$I\times\dg$
along the imbedded $\Sg$ and insert a copy of $M$ so that $\SS^-\sub\hg$. If we
identify $I\times\dg$
with the
complement of the trivial framed
string link then our newly constructed manifold is identified with the
complement of
some string link which lies in a homology $3$-ball precisely when  $M\in\bhg$.
We take this
to
represent $\rho (M)$. It is not hard to see that $\rho\circ\hp =\text{id}$.
\end{proof} 

 \begin{remark}\lbl{rem.prod} Note that $\bhg =\eta_1\i (A)$, where $A$ is the
subset of $\sp
(H)$
consisting of all $h$ such that $p\circ h\circ i$ is an isomorphism. It is not
hard to see that
$A\iso
P\times Q$, where $P, Q$ are subgroups of $\sp (H)$ defined as follows: $h\in P$
if
and only if $h(L)=L$, where $L=\ker p$, and $h\in Q$ if and only if
$h\arrowvert H'=\text{id}$. The
bijection is defined by multiplication $(h_1 ,h_2 )\mapsto h_1\cdot h_2$. This
product decomposition of $A$ lifts to a product decomposition of $\bhg$---see
Corollary \ref{cor.prod}.
\end{remark}
\begin{remark} Since $\rho$ is a retraction, $\rho (\bhg\cap\Gg )$ certainly
contains
$\pgf$, but also contains, for example, the Whitehead string link,
Figure \ref{fig.wh}. 

\begin{figure}[ht!]
\psdraw{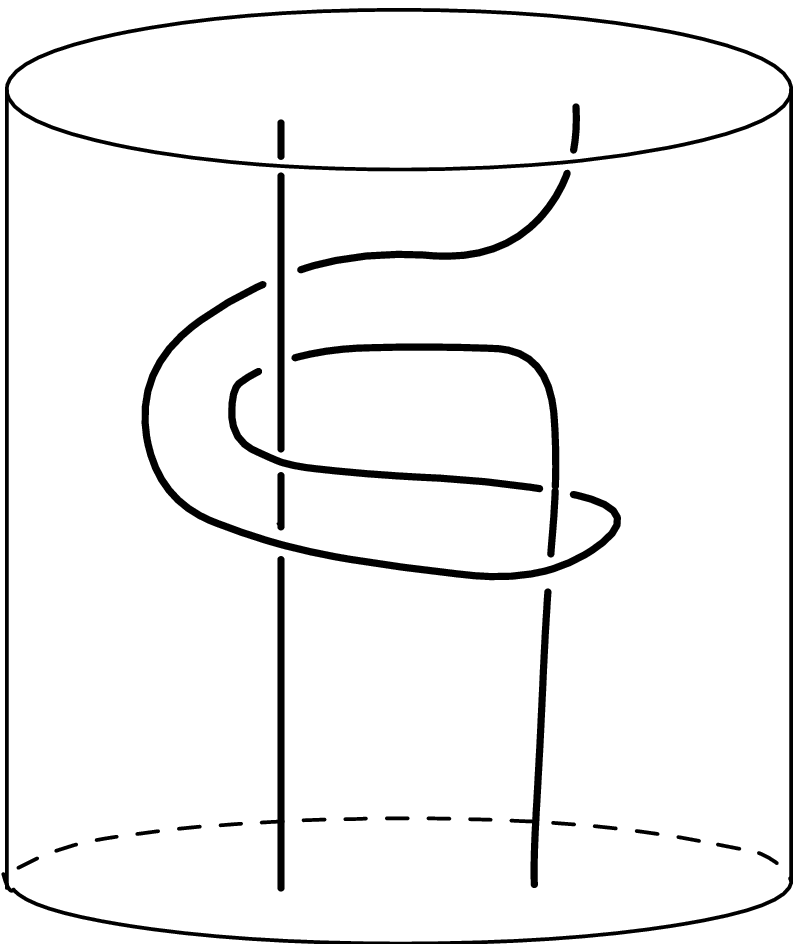}{1.2in}
\nocolon\caption{}\lbl{fig.wh}
\end{figure}

Generally $\s\in\rho (\bhg\cap\Gg )$ if and only if $\s$ is represented by a
framed
string link whose complement is diffeomorphic to $\hg$ (i.e.\ the complement of
the
trivial string link).
\end{remark}

\begin{question} What is $\hp\i (\Gg )$? Clearly $\pgf\sub \hp\i (\Gg )\sub \rho
(\bhg\cap\Gg )$.
\end{question}

\subsection{Boundary homology cylinders}

Let $M$ be a homology cylinder over $\Sg$ and so $\bd M=\SS^+\cup\SS^-$. We can
also write $\SS^+\cup\SS^- =\bd H^+\amalg H^-$, where  $H^{\pm}$ are two copies
of $\hg$ and $\amalg$ denotes boundary connected sum. Let $V$ be the closed
manifold $M\cup (H^+\amalg H^- )$. We will say that $M$ is a {\em boundary
homology cylinder} if $V$ bounds a compact orientable manifold $W$ such that the
inclusions $H^{\pm}\sub W$ are homology equivalences. the boundary homology
cylinders define a subgroup $\bg\sub\Hg$. Note that $\eta_1 (\bg )\sub P$ and so
$\bg\sub\bhg$. In fact,  since $P\cdot A\sub A$, we have $\bg\cdot\bhg\sub\bhg$, 
and so
the right coset space
$\bhg /\bg$ is defined.

\begin{theorem}\lbl{th.bg}
\begin{enumerate}
\item[\rm(a)] $\bg\cap\Gg =\Bg$ the subgroup of $\Gg$ consisting of diffeomorphisms
which
extend to diffeomorphisms of $\hg$.
\item[\rm(b)] $\bg =\rho\i (0)$, where $0$ denotes the trivial string link.
\end{enumerate}
\end{theorem}

\begin{proof}
(a) Suppose $h\in\bg\cap\Gg$. Then $h_* (x_i )\in\ker\{\pi_1 (\Sg )\to\pi_1 (\hg
)\}$ and, therefore, by Dehn's lemma, if $D_i$ is the meridian disk of $\hg$
corresponding to $x_i$, $h(D_i )$ bounds a disk $D'_i\sub\hg$ By standard cut
and
paste techniques we can assume that the $D'_i$ are disjoint. We can now extend
$h$ over
each $D_i$ by mapping it onto $D'_i$. Since the complement of $\cup_i D_i$ is a
$3$-ball, as is the complement of $\cup_i D'_i$, we can extend over $\hg$.

\medskip
(b) Suppose $M\in\bg$. Let $V$ be a homology bordism from $\hg$ to itself,
with
$\bd
V= H^+\cup M\cup H^-$, where $ H^{\pm}$ are two copies of
$\hg$. Let
$W=I\times (I\times D_g )\cup V$, where $V$ is attached to $(1\times\hg )\sub
(1\times (I\times D_g ))$ along $
H^+$. Then $W$ is a homology  bordism between $0\times (I\times D_g )$, the
complement of the trivial string link ({\em triv} in Figure \ref{fig.triv}), and
$C=\text{ complement of }S$, where $S$ is the string link constructed as the
representative of $\rho (M)$. Thus we can fill in $W$ with product strings to
yield a concordance
from the
trivial string link to $S$. See Figure \ref{fig.triv}
\begin{figure}[ht!]
\psdraw{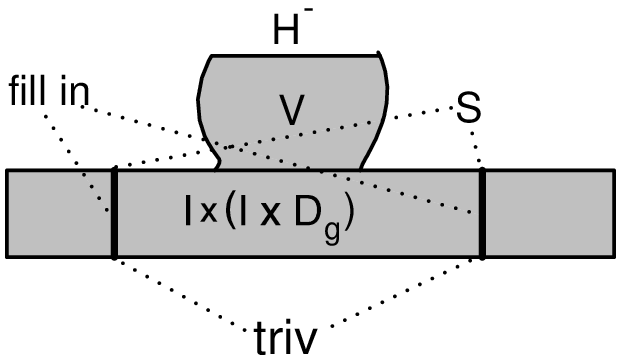}{2in}
\caption{$\rho (\bg )=0$}\lbl{fig.triv}
\end{figure}

Conversely suppose $M\in\bhg$ and $\rho (M)=0$. Let $W$ be the complement of a
concordance between the trivial string link and $S$, the constructed
representative of
$\rho (M)$. If $C$ is the complement of $S$, then we can decompose
$C=I\times\Sg\cup M\cup\hg$. Thus 
$$\bd W=\hg\cup C= H^+\cup M\cup H^- $$
In this way we can see that $W$ is a homology bordism from $\hg$ to itself
which
exhibits $M$ as an element of $\bg$.
\end{proof}

\begin{corollary}\lbl{cor.prod}
The monomorphisms $\bg\sub\bhg$ and $\hp \co\sgf\sub\bhg$ define a bijection
$$\chi \co\bg\times\sgf\to\bhg $$
In particular $\hp$ induces a bijection $\sgf\iso\bhg /\bg$.
\end{corollary}
\begin{proof} $\chi$ is defined by $\chi (M,S)=M\cdot\hp (S)\in\bhg$ (see Remark
\ref{rem.prod}). To see that $\chi$ is one-one we only need note that
$\bg\cap\hp
(\sgf )=1$, which follows from Theorem \ref{th.bg}(b) and that $\rho\circ\hp
=1$. 
To prove that $\chi$ is onto we need:
\begin{lemma}\lbl{lem.hom}
If $M\in\bhg ,S\in\sgf$, then $\rho (M\cdot\hp (S))=\rho (M)\cdot S$.
\end{lemma}
\begin{proof}
The idea is given by the schematic pictures below. The first picture shows the
complement of  $\rho (M\cdot\hp (S))$.  $C$ is the complement of the
string link $S$ and the shaded region is a product $I\times\Sigma_g$. Note that
$M\cdot\hp (S)$ can be constructed by stacking $C$ on $M$ along $D_g\sub\SS^+$.
The second picture gives an alternative view of the first picture which can
then be recognized as the complement of  $\rho (M)\cdot S$. Filling in the
strings completes the proof.

$$
\psdraw{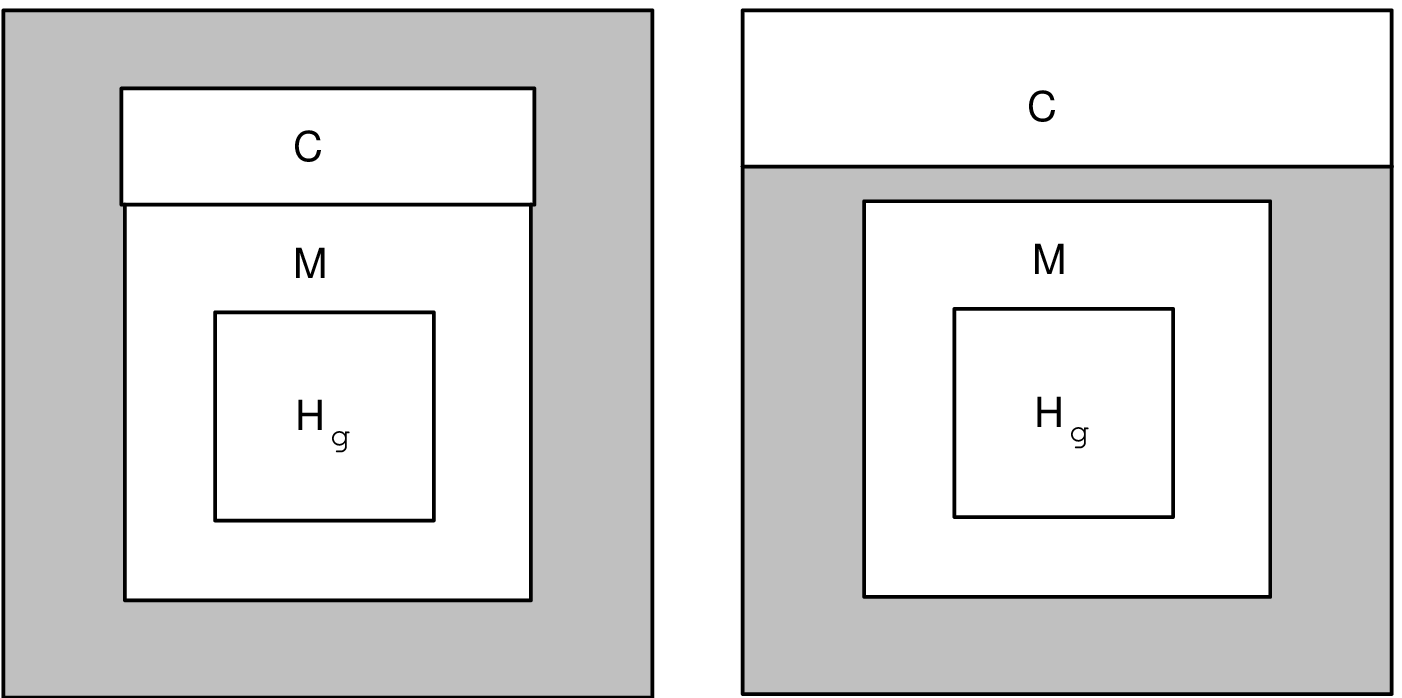}{3.5in}
$$\vglue -20pt
\end{proof}

Now suppose $M\in\bhg$. Then $\rho (M\cdot\hp ((\rho (M))\i )=0$, by Lemma
\ref{lem.hom}, and $\chi (M\cdot\hp ((\rho (M))\i ,\rho (M))=M$.
\end{proof}

\begin{remark} In a recent work of N. Habegger \cite{Ha} a bijection is
constructed between the ``Torelli'' subsemigroup of $\cg$ and the semigroup
of $2g$-strand string links in homology balls with vanishing linking numbers,
which although not multiplicative, induces isomorphisms between the associated
graded groups of the relative weight filtrations.
\end{remark}

The homomorphisms $\eta_k \co\Hg\to\ao (F/F_{k+1})$, induce homomorphisms\nl
$\eta_k^B \co\bg\to\aut (F'/F'_{k+1})$.

\begin{proposition}\lbl{prop.B}
$\eta_k^B$ is onto, for all $k\ge 1$.
 \end{proposition}
 
 \begin{proof}
 If $h\in\aut (F'/F'_{k+1} )$, then we can lift $h$ to an endomorphism $\bar h$
of
$F'$ and we can realize $\bar h$ by an imbedding $\hg\sub\hg$. If $\hg '$
denotes
the
imbedded copy of $\hg$, then define $V=(I\times\hg ')\cup (I\times\hg )$,
attaching
$1\times\hg '$ to $0\times\hg$. $V$ is a homology bordism between $0\times\hg
'$
and $1\times\hg$ and clearly maps to $h$ under $\eta_k^B$.
\end{proof}

\begin{remark}  The restriction of $\eta_k^B$ to $\Gg^B$ is induced by a map
$\eta^B
\co\Gg^B\to\aut F'$ and it is known (see \cite{G}) that $\eta^B$ is onto.
\end{remark}

To see that $\ker\eta_k^B$ is large we consider another imbedding of $\sgf$ into
$\Hg$, $\hat\Psi \co\sgf\to\Hg$, ``dual'' to $\hp$. This is defined in the same
way as
$\hp$---from an imbedding $\dg\sub\Sg$---except that we now choose the imbedding
so that the meridians of the holes in $\dg$ map onto the meridians of the
handles
corresponding to $\{ x_i\}\sub F$. If $F''$ is the free group with basis $\{
x_i\}$, then
we have an injection $\hat\psi \co\A_1 (F''/F''_q )\to\ao (F/F_q )$, defined 
by 
\begin{equation}\lbl{eq.psi}
\begin{split}
\hat\psi ((\l_i ))\cdot x_i &\con\l_i\i   x_i\l_i\mod F_q\\
\hat\psi ((\l_i ))\cdot y_i &\con\l_i y_i \mod F_q
\end{split}
\end{equation}
 There is a commutative diagram
\begin{equation}\lbl{di.psi}
\begin{diagram}
\node{\sgf}\arrow{e,t}{\hat\Psi}\arrow{s,l}{
\mu'_k}\node{\Hg}\arrow{s,r}{\eta_k}\\
\node{\A_1 (F''/F''_{k+1})}\arrow{e,t}{\hat\psi}\node{\ao (F/F_{k+1})}
\end{diagram}
\end{equation}

\begin{theorem}
\begin{enumerate}
\item[\rm(a)] $\hat\Psi (\sgf )\sub\Hg^B$.
\item[\rm(b)] $\eta_k^B\circ\hat\Psi =0$ for all $k$.
\end{enumerate}
\end{theorem}
\begin{proof} 
(a) Write $\hg =D^3\cup H$, where $H$ represents the $g$ disjoint handles. Then
the
imbedding $\dg\sub\Sg$ which defines $\hat\Psi$ can be regarded as the
restriction of
an identification of $D^2$ with the upper hemisphere $D^2_+$ of $D^3$---the
holes
in $D^2$ which are removed to obtain $\dg$ are the disks along which one end of
each
handle of $H$ is attached to $D^3$. Now consider $I\times D^2\sub I\times D^3$
and a
framed string link $S\sub I\times D^2$. We create a manifold $V$ by attaching
$I\times
H$ to $I\times D^3$ by attaching one end of each handle to $S$ and the other to
the
trivial framed string link in $I\times D^2_-$. It is not difficult to see that
this is a
homology bordism from $\hg$ to itself which extends $\hat\Psi (S)$.

\noindent (b) This follows directly from Equation \eqref{eq.psi} and diagram
\eqref{di.psi}, since $\l_i\in F''$.
\end{proof} 
\begin{question}
Describe $\bigcap_k\ker\eta_k^B$. It is known (see \cite{Lu}) that the kernel of
$\eta^B$ is generated by so-called {\em twist automorphisms}, i.e.\ Dehn twists
along
properly imbedded $2$-disks in $\hg$.
\end{question}

\section{A Lagrangian filtration of the mapping class group}

\subsection{Definition of the filtration}
Let $\hg$ and $\Sg$ be as above. We define a new filtration of $\Gg$. Let 
$p\co F\to F'$ 
the
epimorphism induced by the inclusion $\Sg\sub\hg$, i.e.\ $p(x_i )=1,\ p(y_i
)=y_i$.
 We
define $\Lg [k]\sub\Gg$, for $k\ge 1$, by the condition that $h\in\Lg
[k]$ if and only if $p\circ h(x_i )\in F'_{k+1}$ for all $i$. It is not hard to
see that
these are
subgroups. Set $\Lg [\infty
]=\bigcap_k\Lg [k]=\{ h |\ker p\sub\ker p\circ h\}$. Note that $h\in \Lg [k]\
(k\le\infty$)
induces an
automorphism of $F'/F'_{k+1}$ and, therefore, of $H'=H_1 (F')$. This 
defines a
homomorphism $\tau_k \co\Lg [k]\to\auto (H')$. We define $\hat\Lg
[k]=\ker \tau_k$.

We identify some of these groups. Let $L=\ker\{ H\to H'\}$, a 
Lagrangian
subgroup of the symplectic space $H$. Let $\sp (H)$ denote the group 
of
symplectic automorphisms of $H$. Let $ \tau \co\Gg\to\sp (H)$ denote the
standard epimorphism. Let $P(L)\sub\sp (H)$ be the subgroup of all 
$\s$
such that $\s (L)=L$ and $P_0 (L)\sub P(L)$ those $\s$ such that $\s 
|L
=\text{id}_L$. then $\Lg [1]= \tau\i (P(L))$ and $\hat\Lg [1]= \tau\i 
(P_0
(L))$, which was denoted $\bar\Lg^L$ in \cite{GL1}.

\begin{proposition}\lbl{prop.lg}
\begin{itemize}

\item[\rm(a)] $ \tau_k$ is onto for all $k\le\infty$.
\item[\rm(b)] $\hat\Lg [2]$ is the subgroup generated by Dehn 
twists on
simple closed curves representing elements of $L$. This is the 
subgroup denoted
$\Lg^L$ in \cite{GL1}.

\item[\rm(c)] $\Lg [\infty ]=\Gg^B$.
\end{itemize}
\end{proposition}

\begin{proof} (a) Since $\Gg^B\sub\Lg [k]$ for all $k$, this follows
immediately from the fact that $\eta^B$ is onto \cite{G}.

\noi (b) Note that $h\in\hat\Lg [2]$ iff. $h(x_i )x_i\i\in
(F_3\cdot\LL )\cap F_2$ for all $i$, where $\LL =\ker p$, the normal closure of
$\{ x_i\}$. To see this first note that $h\in\Lg [2]$ iff. $h(x_i
)x_i\i\in F_3\cdot\LL$ for every $i$. Secondly note that $h\in\ker\tau_k$ iff.
$h(y_i )\con y_i\mod F_2$ for all $i$---but this is equivalent to
$h(x_i )\con x_i\mod F_2$ for all $i$, since, if a symplectic matrix has the
form
$\left(\begin{array}{cc}A &B\\0&C\end{array}\right)$, it follows that $A=I$ iff.
$C=I$.

Now note that, as a consequence of 
\cite[Cor. 2]{L},  
$h\in\Lg^L$ iff. $h(x_i )x_i\i\in [F,\LL\cdot F_2 ]$ for all $i$. 
Thus, since  $(F_3\cdot\LL
)\cap F_2 =F_3\cdot (F_2\cap\LL )$ and $[F,\LL\cdot F_2 ]=F_3\cdot 
[F,\LL ]$, {\rm (b)} will follow from

\begin{lemma} $F_2\cap\LL =[F,\LL ]$.
\end{lemma}
\begin{proof} By Hopf's theorem the quotient 
$\dfrac{F_2\cap\LL}{[F,\LL ]}=H_2
(F/\LL )$. But $F/\LL\,\iso F'$, which is a free group and so $H_2 
(F/\LL
)=0$.
\end{proof}

\noi (c) If $h\in\Lg [\infty ]$, then $h_* (x_i )\in\LL$, for all $i$. We
can
thus apply the argument in the proof of Theorem \ref{th.bg}(a). 
\end{proof}

Define a function $J_k^L \co\Lg [k]\to\hom (L,\L_{k+1} (H'))\iso 
H'\otimes\L_ {k+1}(H')$ by $\jk (h)\cdot a=p( h(\a ))$, where $\a\in\LL$ is any
lift of
$a\in L$. Note that the symplectic form on $H$ induces an 
identification of
 $H'$ with the dual space of $L$. Furthermore if $h_* \co H\to
H$ is the induced automorphism, then $h_* (L)=L$ and so there is an 
induced
automorphism of $H'$ and of $\L_k (H')$, both of which we also denote 
by
$h_*$.

\begin{proposition}\lbl{prop.J}
\begin{enumerate}
\item[{\rm (a)}] $\jk (\hat\Lg [k])\sub\D_k (H')$.

\item[{\rm (b)}] $\jk (h_1\circ h_2 )=(h_{2*}\otimes 1)\jk
(h_1 )+(1\otimes h_{1\ast})\jk (h_2 )$. Therefore $\jk |\hat\Lg [k]$ 
is a
homomorphism.

\item[{\rm (c)}] $\Lg [k+1]=\ker\jk$.
\end{enumerate}
\end{proposition}

\begin{proof}  (a) If $h\in\Lg [k]$ then $\jk (h)=\sum_i y_i\otimes 
ph_*
(x_i )$. We will abuse notation and allow ourselves to
 denote the induced bases of $L$ and $H'$ by $\{ x_i\}$ and $\{
y_i\}$, respectively. If $h\in\hat\Lg [k]$, then $ph_* (y_i )=y_i$ in
$H'$ and so 
\begin{align*}
\g'_k \jk (h)&=\prod_i [y_i ,ph_* (x_i )]=\prod_i [ph_* (y_i ),ph_*
(x_i )]\\
&=ph_*\prod_i [y_i ,x_i ]=p\prod_i [y_i ,x_i ]=1
\end{align*}
in $F_{k+2}/F_{k+3}$\,, since $h_*\prod_i
[y_i ,x_i ]=\prod_i [y_i ,x_i ]$.

\medskip
\noi (b) Let $\a\in\LL$ be a lift of $a\in L$. Then we can write $h_2 
(\a
)=\a_{1}\a'$ for some $\a_{1}\in\LL$ and $\a'\in F'\sub F$, where the
latter inclusion is some splitting of $p$. Since $h_2\in\Lg [k]$ we 
can
choose $\a'\in F'_{k+1}$. If $a_{1}$ is the homology class of $\a_{1}$ in
$L$, then $a_{1}=h_{2*}(a)$ and so $ph_1 (\a_{1})$ represents $\jk 
(h_1 )
\cdot (h_{2*}a)$. If $\a'$ represents $a'\in F'_{k+1} /F'_{k+2}=\L_{k+1} 
(H')$,
then $a'=\jk (h_2 )\cdot a$. Thus $ph_1 (\a' )$ represents 
$h_{1*}(\jk (h_2
)\cdot a)$.

From these observations we conclude that $\jk (h_1\circ h_2 )\cdot 
a$, which is
represented by $p(h_1 h_2 (\a ))=ph_1 (\a_{1})ph_1 (\a')$ is given by

$$\jk (h_1\circ h_2 )\cdot a=\jk (h_1 )
\cdot h_{2*}a+h_{1*}\jk (h_2
)\cdot a $$
\medskip
\noi (c) This is immediate.
\end{proof}

\subsection{An estimate of $\im\jk$}
We make use of the imbedding $\Phi \co\pgf\to\Gg$.

\begin{theorem}
$\Phi$ induces imbeddings
$$\Phi_k \co(\pgf )_{k+1} /(\pgf )_{k+2}\hookrightarrow\hat\Lg [k]/\hat\Lg
[k+1]\iso\Lg [k]/\Lg [k+1]\hookrightarrow\D_k (H')\quad (k\ge 2)
$$ where $(\pgf )_q$ is the $q$-th lower central series subgroup of
$\pgf$.

The composite imbedding $( \pgf )_ {k+1}/( \pgf )_{k+2}\hookrightarrow\D_k (H')$
is described as follows. If $\b\in (\ \pgf )_ {k+1}$ is
specified by longitude elements $\l_1 ,\dots ,\l_g\in F'_ {k+1}$, then 
$\b$ maps
to $\sum_i y_i\otimes [\l_i ]$, where $y_i$ is a basis of $H'$ and 
$[\l _i
]\in F'_ {k+1}/F'_{k+2}\iso\L_ {k+1}(H')$ denotes the coset of $\l_i$.
\end{theorem}

 From this we conclude that $\rk\Lg [k]/\Lg [k+1]\ge r(g,k)$, where
$r(g,k)=\rk( \pgf )_ {k+1}/( \pgf )_{k+2}$ is explicitly computable
(see the discussion in \cite{L}). The gap between $r(g,k)$ and
$\ker\g_k \co H'\otimes\L_ {k+1}(H')\to\L_{k+2}(H')$ is also explicitly 
computable. For
example if $k=2$ it is $\frac16 (g^3 -g)$ and for $k=3$ it is 
$\frac18 (g^3
-g)(g-2)$---for $k=1$ it is zero. 

\begin{proof} Recall the definition of $\phi$ from Equation \eqref{eq.phi}. 

It is well-known that $\b\in (\pgf )_q$, the $q$-th term of the lower 
central
series, if and only if every $\l_i\in F'_q$---see e.g.\ \cite{F}. We therefore 
have
$\Phi ((\pgf )_ {k+1})\sub\Lg [k]$ and the  induced $\Phi_k$ is injective. That 
$\hat\Lg [k]/\hat\Lg
[k+1]\iso\Lg [k]/\Lg [k+1]$ follows from Proposition \ref{prop.lg}(a).

It 
also
follows directly from the definitions that the composite
$( \pgf )_ {k+1}\to \Lg [k]\to H'\otimes\L_ {k+1}(H')$ is as claimed.
\end{proof}

\section{Lagrangian filtration of $\Hg$}

We can extend the Lagrangian filtration over $\Hg$ in a natural way, using the
$\{\eta_k\}$.  Set $\F^L_k (\Hg )=\{ M\in\Hg \co p\eta_k (M)\arrowvert F''=0\}$,
for
$k\ge
1$, where $p\co F/F_ {k+1}\to F'/F'_ {k+1}$ is the projection and $F''\sub F$
via
$x_i\to x_i$. These are subgroups and obviously $\F^w_k (\Hg )\sub\F^L_k (\Hg ),
\
\F^L_k
(\Hg )\cap\Gg
=\Lg [k]$ and $\Hg^B\sub\F^L_k (\Hg )$ for every $k$. Now $\eta_k$ induces a map
$\eta'_k
\co\F^L_k (\Hg )\to\aut (F'/F'_ {k+1})$ and clearly $\eta'_k\arrowvert\Hg^B
=\eta_k^B$.
Thus $\eta'_k$ is onto by Proposition \ref{prop.B}.

We can also define $\hat\F^L_k (\Hg )=\F^L_k (\Hg )\cap\ker\hat\t$, where
$\hat\t
\co\Hg\to\sp (H)$
is the obvious extension of $\t$.

The map $J_k^L \co\Lg [k]\to H'\otimes L_{k+1} (H')$ extends to a map $\hat
J_k^L \co\F^L_k
(\Hg
)\to H'\otimes L_{k+1} (H')$ and the obvious generalization of Proposition
\ref{prop.J} holds. In
particular
$\F^L_{k+1} (\Hg )=\ker\hat J_k^L$.

We now address the homomorphism $\hat\Phi \co\sgf\to\Hg$. Recall the filtration
$\{\sgf
[k]\}$ from Section \ref{sec.str}.

\begin{theorem}\lbl{th.sL}
$\hat\Phi (\sgf [k])\sub\hat\F^L_k (\Hg )$ and induces isomorphisms

$$
\sgf [k]/\sgf [k+1]\iso\hat\F^L_k (\Hg )/\hat\F^L_{k+1} (\Hg )\iso\F^L_k (\Hg
)/\F^L_{k+1} (\Hg )\iso\D_k (H')
$$
\end{theorem}
\begin{proof}
The first assertion is clear from the definitions. The rest of the theorem
follows from the
observation that the following diagram is commutative
$$
\divide\dgARROWLENGTH by2
\begin{diagram}
\node{\sgf [k]}\arrow[2]{e,t}{\hat\Phi}\arrow{sse,b}{\mu_k}\node[2]{\hat\F^L_{k}
(\Hg )}\arrow{ssw,b}{\hat J_k^L}\\ \\
\node[2]{\D_k (H')}
\end{diagram}
$$
and the fact that $\mu_k$ is onto.
\end{proof}

We now have $\hat\Phi (\sgf [\infty ])\sub\hat\F^L_{\infty} (\Hg )$ and
$\Hg^B\sub\F^L_{\infty} (\Hg )$. We will see that these two subgroups are
independent
and generate $\F^L_{\infty} (\Hg )$.

Recall the bijection $\chi\co\Hg^B\times\sgf\to\bhg$  from Corollary
\ref{cor.prod}.

\begin{theorem}\lbl{th.prod} 
 $$\chi (\Hg^B\times\sgf [k] )=\F^L_k (\Hg ) $$
 for every $2\le k\le\infty$. In particular $\F^L_k (\Hg )$ is generated by the
independent
subgroups $\Hg^B$ and $\hat\Phi (\sgf [k])$.
\end{theorem}
\begin{proof} 
We need to show that $\chi (\Hg^B\times\sgf [k] )=\F^L_k (\Hg )$. First note
that $\chi (
\Hg^B\times
\sgf [k ])\sub\F^L_k (\Hg )$---this follows from Theorem \ref{th.sL}. For
the
ontoness
we can apply the ontoness argument of Corollary \ref{cor.prod} with the extra
fact
that
$\rho (\F^L_k (\Hg ))\sub\sgf [k]$. This follows from the observation that the
following
diagram is commutative
$$
\begin{diagram}
\node{\bhg}\arrow{e,t}{\rho}\arrow{s,l}{\eta_k}\node{\sgf}\arrow{s,r}{\mu_k} \\
\node{\ao (F/F_ {k+1})}\arrow{e,t}{\hat\rho}\node{\A_1 (F'/F'_ {k+1})}
\end{diagram}
$$
where $\hat\rho$ is defined by $\hat\rho (h)=(ph(x_i ))$.
\end{proof}

 Note that $\hat\Phi (\sgf [\infty ])\sub\ker\eta'_k$ and
$\eta'_k\arrowvert\Hg^B$ is
onto, for every $k$.

\ifx\undefined\bysame
	\newcommand{\bysame}{\leavevmode\hbox to3em{\hrulefill}\,}
\fi

\Addresses\recd
\end{document}